\documentclass[10pt]{article}

\usepackage{graphicx}
\usepackage[all]{xy}
\usepackage{color}
\usepackage{amsmath}
\usepackage{amsfonts}
\usepackage{amsthm}
\usepackage{amssymb}
\usepackage{mathrsfs}
\usepackage{tikz}
\usepackage{tikz-cd}
\usepackage{cite}

\usetikzlibrary{matrix}

\newtheorem{thrm}{Theorem}
\newtheorem{lem}[thrm]{Lemma}

\newtheorem{prop}[thrm]{Proposition}
\newtheorem*{thme}{Theorem E1}

\theoremstyle{remark}
\newtheorem{rem}[thrm]{Remark}
\newtheorem{exm}[thrm]{Example}

\theoremstyle{definition}
\newtheorem{defn}[thrm]{Definition}

\numberwithin{thrm}{section} 

\def \Dj{\mbox{\raise0.3ex\hbox{-}\kern-0.4em D}}

\newcommand{\R}{{\mathbb{R}}}
\newcommand{\Z}{{\mathbb{Z}}}

\newcommand{\Q}{{\mathbb{Q}}}

\newcommand{\K}{{\mathbb{K}}}

\newcommand{\cB}{\mathcal{B}}
\newcommand{\cL}{\mathcal{L}}

\newcommand{\la}{{\lambda}}
\newcommand{\al}{{\alpha}}

\newcommand{\im}{\operatorname{im}}

\newcommand{\id}{{\bf{1}}}

\newcommand{\opint}{\operatorname{op-inter}}

\DeclareMathOperator{\pemodop}{{\mathbf{pmod-op}}}
\DeclareMathOperator{\pemod}{{\mathbf{pmod}}}
\DeclareMathOperator{\Diff}{\mathrm{Diff}}
\DeclareMathOperator{\Ham}{\mathrm{Ham}}

\DeclareMathOperator{\powers}{\mathrm{powers}}
\DeclareMathOperator{\Powers}{\mathrm{Powers}}

\begin{document}

\title{Persistence modules with operators \\ in Morse and Floer theory}

\renewcommand{\thefootnote}{\alph{footnote}}

\author{Leonid Polterovich$^{a}$, Egor Shelukhin$^{b}$, Vuka\v sin Stojisavljevi\'c$^{a}$}

\footnotetext[1]{Partially supported by the European Research Council Advanced grant 338809.}
\footnotetext[2]{Partially supported by NSF grant No. DMS-1128155.}

\maketitle

\begin{abstract}
We introduce a new notion of persistence modules endowed with operators. It encapsulates the additional structure on Floer-type persistence modules coming from the intersection product with classes in the ambient (quantum) homology, along with a few other geometric situations. We provide sample applications to the $C^0$-geometry of Morse functions and to Hofer's geometry of Hamiltonian diffeomorphisms, that go beyond spectral invariants and traditional persistent homology.
\end{abstract}

\tableofcontents

\section{Introduction and main results}

Persistent homology is a new field of mathematics that originates in data analysis. It can be considered a book-keeping device for the information about the topology of sub-level sets in Morse theory (and its generalizations, including Floer theory) that is stable under $C^0$-perturbations of the Morse functions (or the functionals of Floer theory). In the case of Hamiltonian Floer theory, the role of the $C^0$-metric is played by Hofer's metric on the group of Hamiltonian diffeomorphisms. The invariants of persistent homology can be described as a collection of intervals in the real line, called a barcode. In a recent paper \cite{PolShe} (extended in \cite{Zhang}) it was observed that the arithmetical properties of barcodes are pertinent to such questions on Hofer's geometry as the study of the minimal Hofer norm of a perturbation of a given Hamiltonian diffeomorphism necessary to make it autonomous, or more generally - admit a root of order $k \geq 2.$ In the current paper, we introduce and discuss the notion of persistence modules with operators, that allows us to use operators of intersection with cycles in the ambient (quantum) homology to further control the multiplicities of bars in the barcode. This provides new results on Hofer's geometry, and can be shown to provide strictly new information, as compared with traditional persistent homology (including spectral invariants), about the $C^0$-geometry of Morse functions.

\subsection{The Arnol'd conjecture}

As explained in Arnol'd's \cite{MMCM}, an important invariant formulation of the equations of motion of classical mechanics involves a manifold $M,$ {\em the phase space}, a closed non-degenerate two-form $\omega$ on $M,$ {\em the symplectic form}, and a smooth, possibly time-dependent, {\em Hamiltonian} function $H:[0,1] \times M \rightarrow \mathbb{R}$ on $M,$ the total energy function of the system.  The dynamics on the symplectic manifold $(M,\omega),$ which we assume to be closed throughout this paper, is then described by the \textit{Hamiltonian flow} of $H$, \[\{\phi^t_H: M \to M\}_{t\in [0,1]},\] obtained by integrating the time-dependent vector field $X^t_{H}$ on $M$, given by setting $H_t(-) = H(t,-),$ and $$\omega(X^t_{H},\cdot)=-dH_t(\cdot).$$

The time-one map $\phi=\phi_H^1$ of this flow is called a \textit{Hamiltonian diffeomorphism} and it is clear that 1-periodic orbits of the flow $\{\phi^t_H\}_{t\in [0,1]}$ correspond to the fixed points of $\phi$. Hamiltonian diffeomorphisms form a group, that we denote $\Ham(M,\omega).$

In the 1960's Arnol'd has proposed a famous conjecture \cite{ArnoldConjecture2,ArnoldConjecture1} that, essentially, the number of fixed points of $\phi \in \Ham(M,\omega)$ should satisfy the same lower bounds as the number of critical points of a smooth function $f$ on $M.$ The most common interpretation of this conjecture states that under the nondegeneracy\footnote{Meaning that $\det (\id-D\phi_1^H(x))\neq 0$ for every fixed point $x$ of $\phi^1_H$, i.e. the graph of $\phi^1_H$ intersects the diagonal in $M \times M$ transversely.} assumption on $H$, the number of fixed points of $\phi^1_H$ is bounded from below by the sum of the rational Betti numbers of $M$. This conjecture has been a major driving force for the development of the field of symplectic topology. It was first proven in dimension $2$ by Eliashberg \cite{Eliashberg2d}, for tori of arbitrary dimension by Conley and Zehnder \cite{ConleyZehnder}, and on complex projective spaces by Fortune and Weinstein \cite{Fortune,FortuneWeinstein}. The decisive breakthrough on this question was achieved by Floer \cite{Floer1,Floer2,Floer3}, who combined the variational methods of Conley-Zehnder and Gromov's then-recent discovery of the theory of pseudoholomorphic curves in symplectic manifolds \cite{GromovPseudohol}, to construct a homology theory on the loop space of $M$ which parallels the more classical Morse homology (in turn originating in Witten's interpretation of Morse theory \cite{Witten}; cf. \cite{SchwarzMorse}).

Let us briefly describe Floer's work in the simplest setting (cf. \cite{AudinDamian}). Assuming that $M$ is symplectically aspherical, that is $\omega|_{\pi_2(M)}=0, \;c_1(M,\omega)|_{\pi_2(M)} = 0,$ one can define the \textit{action functional} $\mathcal{A}_H$ on the space $\mathcal{L}_c M$ of contractible loops on $M$ by setting
$$\mathcal{A}_H(z)=\int_{0}^{1}H_t(z(t))~ dt - \int_{D^2} {\overline{z}}^* \omega, $$
where $z:[0,1]\rightarrow M$, $z(0)=z(1)$ and $\overline{z}:D^2\rightarrow M$, $\overline{z}(e^{2\pi i t})=z(t)$. Indeed by the asphericity assumption this value depends only on $z,$ and not on $\overline{z}.$

The periodic orbits of the Hamiltonian flow coincide with the critical points of the action functional $\mathcal{A}_H$, which serves as a Morse function in the construction of Floer homology. Since critical points give generators of Morse chain complexes, if $H$ is non-degenerate, the Floer chain complex $CF_*(H)$ will be generated by contractible periodic orbits of $H$. Floer homology will be isomorphic to singular homology of $M$ with rational coefficients, and hence there must be at least
$$\dim CF_*(H)\geq \dim HF_*(H)=\sum_k \dim H_k(M, \mathbb{Q}),$$
such periodic orbits in $M.$ This solves the rational homological version of the Arnol'd conjecture. The detailed construction of Floer homology is rather involved and it has been developed in increasing generality over the years by the combined work of many people (cf. \cite{FukayaOno1,FukayaOno2,PSS,LiuTian,Ruan,PardonHam}), in particular proving the above statement for general $M$ (without assumptions on $\pi_2(M)$). However, various other interpretations of the conjecture are still open in general, with only partial results currently achieved (see e.g. \cite{Floer3,RudyakOprea,FloerCuplength,HoferLusternik,RizellGolovko,OnoPajitnov,Barraud}).

\subsection{Persistent homology}

In both Morse and Floer theory, the differential is given by counting certain trajectories of the negative gradient vector field connecting pairs of critical points. Since the Morse function, or the action functional in the Floer case, decreases along these trajectories, for each $s \in \R$ the generators whose critical values are $<s,$ form a subcomplex. In the case of a Morse function $f$ on a closed manifold $X$, the homology $V^s_*(f)$ of this subcomplex with coefficients in a field $\K$ is isomorphic to the homology $H_*(\{f<s\},\K)$ of the sublevel set $\{f<s\}$ with coefficients in $\K$ (which is a vector space of finite dimension!). Inclusions $\{f<s\} \subset \{f<t\}$ for $s \leq t$ yield maps $\pi_{s,t}: V^s_*(f) \to V^t_*(f)$ such that $\pi_{s,t} \circ \pi_{r,s} = \pi_{r,t}$ for $r \leq s \leq t.$ This family of vector spaces parametrized by a real parameter forms an algebraic structure called a {\em persistence module} which was introduced and studied extensively since the early 2000s in the data analysis community (see e.g. \cite{CEH-stability,CdSGO-structure,ELZ-simplification, BauLes,Carlsson,CZCG,CrawBo,Ghrist,CarlZom}). Quite recently, persistence modules found applications in symplectic topology, see \cite{Team,PolShe,UsherZhang,Zhang,Fraser,Stevenson}, with preludes in \cite{UsherBD1,UsherBD2,FOOO-polydiscs,Barannikov}.
\\

In the case of Floer theory on a symplectically aspherical manifold, the exact same procedure applies. In the so-called monotone case, when $\omega|_{\pi_2(M)} = \kappa \cdot c_1(TM)|_{\pi_2(M)}$ for a constant $\kappa > 0,$ one can define a persistence module by looking at $V^s_m(H)$ for a fixed degree $m,$ which is the approach we choose later in the paper. The structure theorem of persistence modules allows us to associate to this situation a barcode: a multi-set $\mathcal{I}$ of intervals of the form $(a,b]$ or $(a,\infty),$ such that
$$V^t_m(H)\cong \bigoplus_{I\in \mathcal{I}} Q^t(I),$$
where $Q^t(I)$ is a persistence modules equal to $\mathbb{K}$ for $t\in I$ and 0 otherwise. The number of the infinite intervals in the barcode is equal to the dimension of the ambient (quantum) homology in a given degree, while the left ends of these intervals correspond to the well know \textit{spectral invariants}. Spectral invariant can be defined in complete generality (without assumptions on $\pi_2(M)$), and have been used extensively in symplectic topology in the last few decades starting with the foundational papers \cite{Oh-construction,Schwarz:action-spectrum,Viterbo-specGF} (see \cite{LeclercqZapolsky,Usher-spec} for recent developments and numerous further references), before persistence modules entered the field. In terms of barcodes, in complete generality one can only expect bars with endpoints in $\R/\mathcal{P}(\omega),$ where $\mathcal{P}(\omega) = \im \,(\int \omega: \pi_2(M) \to \R)$ is the period group of $\omega.$ In \cite{UsherZhang} such a persistence module is constructed, with, remarkably, the lengths of intervals being well-defined.
\\

Denote the $C^0$-distance between two smooth functions $f$ and $g$ on a compact manifold $X$ by
$$|f-g|_{C^0}=\max_{x \in X} |f(x)-g(x)|,$$
and Hofer's $L^{1,\infty}$-distance between Hamiltonians $F_t$ and $G_t$ on $M$, $t\in [0,1]$ by
$$\mathcal{E}(F-G)=\int_{0}^{1} \bigg( \max_{x \in M} (F_t(x)-G_t(x)) - \min_{x \in M} (F_t(x)-G_t(x)) \bigg) ~dt .$$ The Hofer's pseudo-metric on the universal cover $\widetilde{\Ham}(M,\omega)$ of $\Ham(M,\omega)$  is defined as \[\widetilde{d}(\widetilde{f},\widetilde{g}) = \inf \mathcal{E}(F-G),\] where the infimum runs over all the $F,G$ such that $[\{\phi^t_F\}] = \widetilde{f},$  $[\{\phi^t_G\}] = \widetilde{g}$ in  $\widetilde{\Ham}(M,\omega).$ Similarly, Hofer's metric \cite{HoferMetric,Lalonde-McDuff-Energy} (cf. \cite{P-book,HZ-book}) on $\Ham(M,\omega)$ is  \[d(f,g) = \inf \mathcal{E}(F-G),\] where the infimum runs over all the $F,G$ such that $\phi^1_F = f,$  $\phi^1_G = g$.

One crucial feature of the barcodes is that if the $C^0$-distance between Morse functions, or the Hofer's $L^{1,\infty}$-distance between two Hamiltonians is at most $c,$ then one barcode can be obtained from the other, by a procedure that allows to move each endpoint of a bar by distance at most $c$ (this allows erasing intervals of length $\leq 2c,$ as well as creating such intervals). This follows from the celebrated {\em isometry theorem} of persistence modules. Such a procedure is encoded by the formal notion of a $c$-matching of barcodes, whose algebraic counterpart is the notion of a $c$-interleaving. More precisely, we say that two persistence module morphisms $f_t:V^t \rightarrow W^{t+c}$ and $g_t:W^t \rightarrow V^{t+c}$ induce a $c$-interleaving between persistence modules $V$ and $W$ if $ g_{t+c} \circ f_t = \pi^V_{t,t+2c} $ and $ f_{t+c} \circ g_t = \pi^W_{t,t+2c} $ for every $t\in \mathbb{R}$. The main point is that Morse and Floer continuation maps with respect to linear interpolations of functions $f$ and $g$ or Hamiltonians $F$ and $G$ (or small perturbations thereof) yield the required $c$-interleavings, where $c=|f-g|_{C^0}$ or $c=\mathcal{E}(F-G)$. This allows us to bound the $C^0$-distance between Morse functions as well as Hofer's distance between Hamiltonians (and consequently Hamiltonian diffeomorphisms) from below by the minimal $c$ needed to match the corresponding barcodes. 
\\

In this paper we primarily investigate an additional structure on Morse and Floer persistence modules coming from the ambient homology. Our main observation is that the ambient homology acts on the persistence module by intersecting cycles in the sublevel sets of functions (and a similar picture holds in the Floer case). We consider this action as a particular case of the notion of a persistence module with an operator. Namely, we consider pairs $(V,A)$ where $A:V^t\rightarrow V^{t+c_A}$ is a persistence module morphism as main objects of interest and define morphisms between these objects to be usual persistence module morphisms which commute with the corresponding operators. We may now define operator interleaving as an interleaving in this new category, i.e. an interleaving which commutes with the operators. The fact that $(V,A)$ and $(W,B)$ are $c$-operator interleaved will immediately imply that $\im A$ and $\im B$ (as well as $\ker A$ and $\ker B$) are $c$-interleaved (see Section~\ref{Sect:pmod-operator} for a discussion of persitence modules with operators). 

In the Morse and the Floer case, fixing a (quantum) homology class $a$, we obtain an operator $a*$ induced by intersection (or quantum) product. Continuation maps commute with this operator, hence constitute morphisms of persistence modules with operators and induce operator interleavings. Finally, they provide both $\im(a*)$ and $\ker(a*)$ for two functions $f$ and $g$ or two Hamiltonians $F$ and $G$, with $c$-interleavings, for $c=|f-g|_{C^0}$ or $c=\mathcal{E}(F-G)$ respectively. This means that we may bound these values from below by using barcodes associated to $\im(a*)$ or $\ker(a*)$. Following this line of reasoning, we show that there exists a pair $f,g$ of Morse function on a manifold (even of dimension $2$) such that all their spectral invariants, as well as their barcodes coincide, and yet the corresponding $\im(a*)$ modules are at a positive (computable) interleaving distance $c$. We conclude that the two functions must be at $C^0$-distance at least $c$ (see Section~\ref{Sect: Morse} for an example).
\\

Finally, we present an application to Hofer's geometry, by proving new cases of the conjecture that on any closed symplectic manifold and for any integer $k \geq 2,$ there exist Hamiltonian diffeomorphisms which are arbitrarily far away, in Hofer's metric, from having a root of order $k.$ First results of this kind were obtained in \cite{PolShe}, and were then extended to certain other cases in \cite{Zhang} (for $k$ sufficiently large). In our situation, the multiplication with classes in ambient homology allows to adjust multiplicities of certain long bars, the number theoretic properties of which are crucial to the argument, and allow to reduce the $k$ for which the result holds, yielding Theorem~\ref{main} (see Section~\ref{Sect:Hofer}).

\subsection{Hofer's distance to $k$-th powers}\label{Sect:Hofer}

In this section we shall assume that the ground field $\K$ has characteristic $\mathrm{char}(\K) \neq p,$ contains all $p$-th roots of unity\footnote{That is the polynomial $x^p-1 \in \K[x],$ which is separable by the assumption $\mathrm{char}(\K) \neq p,$ splits over $\K.$}, and fixing a primitive $p$-th root of unity $\zeta_p,$ the equation $x^p -(\zeta_p)^q = 0,$ for each integer $q$ coprime to $p,$ has no solutions in $\K.$ An example of such a field is the splitting field $\Q_p$ over $\Q$ of $x^p-1 \in \Q[x].$\\

Let $(M,\omega)$ be a closed symplectic manifold,  and put $\Powers_k(M,\omega) \subset \Ham(M,\omega),$ where $k$ is an integer, for the set of all diffeomorphisms in $\Ham(M,\omega),$ admitting a root of order $k$ (in the same group). The following result was proven in \cite{PolShe}.

\begin{thrm}\label{Thm:old}
Let $(\Sigma,\sigma)$ be a closed Riemann surface of genus at least $4,$ endowed with an area form, and let $(N,\omega_N)$ be either a point, or a closed symplectically aspherical symplectic manifold. Then for each $k \in \Z_{\geq 2}$ there exists a sequence $\phi_j \in \Ham(\Sigma,\sigma),$ such that \[d(\phi_j \times \id_N, \Powers_k(\Sigma \times N)) \xrightarrow{j \to \infty} \infty.\]
\end{thrm}

Now, assume that $N$ is a monotone symplectic manifold, fix a prime number $p$ and denote by $\Lambda_{\K}$ the field of Laurent power series in variable $q^{-1}$ with coefficients in $\K$,
$$ \Lambda_{\K}= \bigg\{ \sum\limits_{n\in \mathbb{Z}}a_nq^n ~ \bigg\vert ~ a_n \in \K,~(\exists n_0 \in \mathbb{N})~a_n=0 ~\text{for}~n\geq n_0 \bigg\}. $$
The quantum homology of $N$ with $\K$ coefficients is the vector space $H_*(N,\K) \otimes_{\K} \Lambda_{\K}$ over $\Lambda_{\K}$ which we denote by $QH(N)$.
Assuming that $\deg q=2c_N$ where $c_N$ is the minimal Chern number of $N$, $QH(N)$ has a natural $\mathbb{Z}$-grading, that is we can define $QH_r(N)$ for $r\in \mathbb{Z}$, which will be vector spaces over the base field $\K$. We also have that $QH_{r+2c_N}(N) \cong QH_{r}(N)$ for every $r\in \mathbb{Z}$, where the isomorphism is given by multiplication by $q$. Let $e \in QH(N)$ be a homogeneous element and define a map
$$e* : QH(N) \rightarrow QH(N),~ (e*)a=e*a,$$
where $*$ denotes quantum product. This map is a linear morphism between vector spaces over $\Lambda_{\K}$ which restricts to a linear morphism between vector spaces over $\K$ after fixing the grading:
$$e*:QH_r(N) \rightarrow QH_{r-2n+\deg e}(N),~ \text{for}~r\in \mathbb{Z}.$$
Now $E:=e*(QH(N))\subset QH(N)$ is a vector space over $\Lambda_{\K}$ and
$$E_r:=e*(QH_r(N))\subset QH_{r-2n+\deg e}(N),$$
are vector spaces over $\K$ which satisfy $E_r \cong E_{r+2c_N}$, the isomorphism being induced by multiplication by $q$. These spaces give us $2c_N$ Betti numbers associated to a homogeneous element $e\in QH(N)$, which we define as:
$$b_r(e)=\dim_{\K} E_r, ~r=0, \ldots , 2c_N-1.$$
Now, we can state the result regarding Hofer's geometry. Denote by $\powers_p(M)$ the supremum of the Hofer distance to $p$-th powers in $\Ham(M)$. That is for each $\phi \in \Ham(M)$ define $d(\phi,\Powers_p(M)) = \displaystyle\inf_{\theta \in \Powers_p(M)} d(\phi,\theta),$ and define $\powers_p(M) := \sup_{\phi \in \Ham(M)} d(\phi,\Powers_p(M)).$

{\thrm If there exists $e\in QH(N)$ such that $p \nmid b_r(e)$ for some $r\in \{0, \ldots ,2c_N-1 \}$ then
$$\powers_p(\Sigma \times N)=+ \infty .$$
\label{main}}

To prove this result we describe the Floer theoretical setup that fits into our algebraic framework of equivariant persistence modules with operators, and then make a concrete computation in the case of the egg-beater flow which yields the result.

{\exm \label{exa:main-thm} Taking $N$ to be any monotone or symplectically aspherical manifold and $e\in QH(N)$ any class we have $p\nmid b_r(e)$ for large enough $p$. This means that for large enough $p$
$$\powers_p(\Sigma \times N)=+ \infty.$$
Since autonomous Hamiltonian diffeomorphisms have $p-th$ roots for all $p$, we in particular have that Hofer's distance to autonomous flows in $\Ham(\Sigma \times N)$ is unbounded.
}
{\exm \label{exa2:main-thm} Let $N$ be connected, $\dim N=2n$ and assume $c_N\geq n+1$. We now have that $b_0([N])=b_0(N)=1$, where $[N]$ is the fundamental class and $b_0(N)$ classical Betti number, and hence
$$\powers_p(\Sigma \times N)=+ \infty, ~\text{for all}~p.$$
This is for example the case for $N=\mathbb{C}P^n$. Connected symplectically aspherical $N$ fall in this class of manifolds, with $c_N = +\infty.$}
{\exm Let $N=S^2 \times S^2$ and denote by $P=[pt]$, $A=[S^2\times pt]$, $B=[pt \times S^2]$ and by $[N]$ the fundamental class. These four classes form a basis of $QH(N)$ over $\Lambda_{\K}$ and multiplication is completely described by the relations
$$A*B=P,~A^2=B^2=q^{-1}[N].$$
We calculate
$$(A+B)*A=(A+B)*B=P+q^{-1}[N]\in QH_0(N),$$
as well as
$$(A+B)*[N]=A+B \in QH_2(N),~(A+B)*P=q^{-1}(A+B)\in QH_{-2}(N),$$
and hence $b_0(A+B)=b_2(A+B)=1$. This implies that
$$\powers_p(\Sigma \times N)=+ \infty, ~\text{for all}~p.$$
Note that in this example it is crucial that $A+B$ is not invertible. Otherwise, multiplication would be an isomorphism of $QH(N)$ and all the Betti numbers would be equal to 2, so we would have to assume $p\geq 3$.}
{\rem A different extension of \cite[Theorem 1.3]{PolShe}, using different methods, was obtained recently by Zhang in \cite{Zhang}. The result refers to a more general manifold, namely the product $\Sigma \times N$, where $N$ is any symplectic manifold (not necessarily monotone or aspherical) and gives a condition on $p$ in terms of quantum Betti numbers for $\powers_p(\Sigma \times N)$ to be infinite. The $k$-th quantum Betti number is defined as
$$qb_k(N)=\sum\limits_{s\in \mathbb{Z}} b_{k+2c_N\cdot s}(N),$$
where $b_i(N)$ are classical Betti numbers. The main theorem of \cite{Zhang} states that if
$$p \nmid qb_p(N)+2qb_0(N)+qb_{-p}(N),$$
then
$$\powers_p(\Sigma \times N)=+ \infty.$$
One immediately sees that when $N$ is monotone, $qb_k=b_k([N])$, thus in this case our theorem implies Zhang's result. The above examples of $N=\mathbb{C}P^n$ and $N=S^2\times S^2$, show that our criterion sometimes gives a strictly better answer, since the criterion from \cite{Zhang}, fails when $p=2$.

\section{Persistence modules}

\subsection{Basics}\label{Sect:pmod}

We recall briefly the category $\pemod$ of persistence modules that we work with, together with their relevant properties. For detailed treatment of these topics see \cite{BauLes,Carlsson,CZCG,CrawBo,Ghrist,CarlZom}. 
\\
\\
Let $\mathbb{K}$ be a field. A persistence module over $\mathbb{K}$ is a pair $(V,\pi)$ where, $\{V^t\}_{t\in \mathbb{R}}$ is a family of finite dimensional vector spaces over $\mathbb{K}$ and $\pi_{s,t}:V^s \rightarrow V^t$ for $s<t,~s,t\in \mathbb{R}$ is a family of linear maps, called {\it structure maps}, which satisfy:
\begin{enumerate}

  \item[1)] $V^t=0$ for $t \ll 0$ and $\pi_{s,t}$ are isomorphisms for all $s,t$ sufficiently large;
  \item[2)] $\pi_{t,r}\circ \pi_{s,t}=\pi_{s,r}$ for all $s<t<r$;
  \item[3)] For every $r\in \mathbb{R}$ there exists $\varepsilon >0$ such that $\pi_{s,t}$ are isomorphisms for all $r-\varepsilon <s <t \leq r$;
  \item[4)] For all but a finite number of points $r\in \mathbb{R}$, there is a neighbourhood $U\ni r$ such that $\pi_{s,t}$ are isomorphisms for all $s<t$ with $s,t \in U$.
\end{enumerate}
The set of the exceptional points in 4), i.e. the set of all points $r\in \mathbb{R}$ for which there does not exist a neighbourhood $U\ni r$ such that $\pi_{s,t}$ are isomorphisms for all $s,t \in U$, is called {\it the spectrum} of the persistence module $(V,\pi)$ and is denoted by $\mathcal{S}(V)$. One easily checks that for two consecutive points $a<b$ of the spectrum and $a<s<t \leq b$, $\pi_{s,t}$ is an isomorphism. This means that $V^t$ only changes when $t$ "passes through points in the spectrum".

We define a morphism between two persistence modules $A:(V,\pi) \rightarrow (V',\pi')$ as a family of linear maps $A_t:V^t \rightarrow (V')^t$ for every $t\in \mathbb{R}$ which satisfies $$A_t \pi_{s,t}=\pi_{s,t}'A_s ~ ~ \text{for} ~ s<t.$$ Note that the kernel $\ker A$ and an image $\im A$ are naturally persistence modules whose families of vector spaces are $\{\ker A_t \subset V^t\}_{t \in \R},$ $\{\im A_t \subset (V')^t\}_{t \in \R},$ since the structure maps $\pi_{s,t}$ restrict to these systems of subspaces. In fact, it is not difficult to prove that $\pemod$ forms an abelian category, with the direct sum of two persistence modules $(V,\pi)$ and $(V',\pi')$ given by
 $$(V,\pi)\oplus(V',\pi')=(V\oplus V', \pi \oplus \pi').$$

{\exm Let $X$ be a closed manifold and $f$ a Morse function on $X$. For $t\in \mathbb{R}$ define $V^t(f)=H_*(\{ f<t \},\mathbb{K})$ to be homology of sublevel sets of $f$ with coefficients in a field $\mathbb{K},$ and let $\pi_{s,t}:V^s(f) \rightarrow V^t(f)$ be the maps induced by inclusions of sublevel sets. One readily checks that $(V(f),\pi)$ is a persistence module. The spectrum of $V(f)$ consists of critical values of $f$. Similarly fixing a degree $r \in \Z,$ one obtains a persistence module $V^t_r(f)=H_r(\{ f<t \},\mathbb{K}).$ It is easy to see that the spectrum of $V_r(f)$ is contained in the set of critical values of $f$ of critical points of index $r$ or $r+1.$ Finally $V(f) = \oplus V_r(f).$ Hence $V(f)$ has the structure of a persistence module of $\Z$-graded vector spaces.\label{BasicExample}}

\medskip

An important object in our story is the barcode associated to a persistence module. It arises from the structure theorem for persistence modules, which we now recall. Let $I$ be an interval of the form $(a,b]$ or $(a,+\infty)$, $a,b\in \mathbb{R}$ and denote by $Q(I)=(Q(I),\pi)$ the persistence module which satisfies $Q^t(I)=\mathbb{K}$ for $t\in I$ and $Q^t(I)=0$ otherwise and $\pi_{s,t}=id$ for $s,t\in I$ and $\pi_{s,t}=0$ otherwise.

 {\thrm[The structure theorem for persistence modules] For every persistence module $V$ there is a unique collection of pairwise distinct intervals $I_1,\ldots , I_N$ of the form $(a_i,b_i]$ or $(a_i,+\infty)$ for $a_i,b_i \in \mathcal{S}(V)$ along with the multiplicities $m_1,\ldots,m_N$ such that
 $$V\cong \bigoplus\limits_{i=1}^N(Q(I_i))^{m_i}.$$
 \label{StructureTheorem}}
The multi-set which contains $m_i$ copies of each $I_i$ appearing in the structure theorem is called {\it the barcode} associated to $V$ and is denoted by $\mathcal{B}(V)$. Intervals $I_i$ are called {\it bars}.

\medskip

\begin{rem}
One feature of the Example~\ref{BasicExample} is the existence of additional structure that comes from identifying  $V^\infty := \displaystyle\varinjlim V^t$ with $H_*(X,\mathbb{K}).$ Put $\Psi: V^\infty \to H_*(X,\mathbb{K})$ for the natural isomorphism. Given $a \in H_*(X,\mathbb{K})$ with $a \neq 0,$ we can produce the number $c(a,f) := \inf\{t \in \R \,|\, \Psi^{-1}(a) \in \im (V^t \to V^\infty)\}.$ This number is called a {\em spectral invariant}, and has many remarkable properties. One can prove that for each $a \neq 0,$ $c(a,f)$ is a starting point of an infinite bar in the barcode of $V(f),$ and each such starting point can be obtained in this way.
\end{rem}

\medskip

For an interval $I=(a,b]$ or $I=(a,+\infty)$, let $I^{-c}=(a-c,b+c]$ or $I^{-c}=(a-c,+\infty)$, and similarly $I^c=(a+c,b-c]$ or $I^c=(a+c,+\infty),$ when $b-a>2c$. We say that barcodes $\cB_1$ and $\cB_2$ admit a $\delta$-matching if it is possible to delete some of the bars of length $\leq 2\delta$ from $\cB_1$ and $\cB_2$ (and thus obtain $\bar{\cB}_1$ and $\bar{\cB}_2$) such that there exists a bijection $\mu:\bar{\cB}_1 \rightarrow \bar{\cB}_2$ which satisfies
$$\mu(I)=J \Rightarrow I\subset J^{-\delta},~J\subset I^{-\delta}.$$
We define {\it the bottleneck distance} $d_{bottle}(\cB_1,\cB_2)$ between barcodes $\cB_1,\cB_2$ as infimum over $\delta >0$ such that there exists a $\delta$-matching between them. One readily checks that the following lemma holds.
{\lem Let $V_1,\ldots, V_l$ and $W_1,\ldots, W_l$ be persistence modules. Then
$$\mathcal{S}(\bigoplus_{r=1}^l V_r)=\bigcup\limits_{r=1}^l\mathcal{S}(V_r),~~~ \mathcal{B}(\bigoplus_{r=1}^l V_r)=\sum\limits_{r=1}^l\mathcal{B}(V_r),$$
and
$$d_{bottle}\bigg( \mathcal{B}(\bigoplus_{r=1}^l V_r),\mathcal{B}(\bigoplus_{r=1}^l W_r) \bigg)\leq \max_r d_{bottle}(\mathcal{B}(V_r),\mathcal{B}(W_r)).$$
Here $\Sigma$ denotes multiset sum, that is union of elements, adding up multiplicities. \label{sumpersistence}}
\\
\\
For a persistence module $V=(V,\pi)$ denote by $V[\delta]=(V[\delta],\pi[\delta])$ a shifted persistence module given by $V[\delta]^t=V^{t+\delta},\pi_{s,t}=\pi_{s+\delta,t+\delta}$ and by $sh(\delta)_{V}:V\rightarrow V[\delta]$ a canonical shift morphism given by $(sh(\delta)_{V})_t=\pi_{t,t+\delta}:V^t \rightarrow V^{t+\delta}$. Note also that a morphism $f:V\rightarrow W$ induces a morphism of $f[\delta]:V[\delta] \rightarrow W[\delta]$. We say that a pair of morphisms $f:V\rightarrow W[\delta]$ and $g:W\rightarrow V[\delta]$ is a $\delta$-interleaving between $V$ and $W$ if
$$g[\delta]\circ f = sh(2\delta)_{V}\text{ and } f[\delta]\circ g=sh(2\delta)_{W}.$$
Now we can define {\it the interleaving distance } $d_{inter}(V,W)$ between $V$ and $W$ as infimum over all $\delta>0$ such that $V$ and $W$ admit a $\delta$-interleaving. The isometry theorem for persistence modules states that $d_{inter}(V,W)=d_{bottle}(\mathcal{B}(V),\mathcal{B}(W))$ (see \cite{BauLes}).

\subsection{K\"unneth formula for persistence modules}\label{Sect:Kunneth}

As we mentioned before, $\pemod$ is an abelian category, and we wish to define a monoidal structure $\otimes$ and its derived functors in this category in a similar fashion to the situation which we have for $\mathbb{Z}$ modules (similar constructions, yet with different aims and applications, appeared in \cite{CarlZom2,CarlZom3,CurryPhD,VongPhD,VongCarl}).
Let $(V^s,\pi^V)$ and $(W^t,\pi^W)$ be two persistence modules and define vector spaces
$$X^r=\bigoplus_{t+s=r} V^s \otimes W^t,~\text{and}~Y^r\subset X^r~\text{for every}~r\in \mathbb{R},$$
given by
$$Y^r= \bigg\langle \bigg\{ (\pi^V_{\alpha,s_1} v_\alpha)\otimes (\pi^W_{\beta,t_1} w_\beta) - (\pi^V_{\alpha,s_2} v_\alpha)\otimes (\pi^W_{\beta,t_2} w_\beta) \bigg\} \bigg\rangle,$$
where $\langle S \rangle$ stands for vector space over $\mathbb{K}$ generated by the set $S$ and indices $s_1,s_2,t_1,t_2,\alpha$ and $\beta$ satisfy $s_1+t_1=s_2+t_2=r, \alpha \leq \min \{s_1,s_2 \},~\beta \leq \min \{t_1,t_2 \}$.
We may now define $(V\otimes W)^r=X^r / Y^r$ and maps $\pi^V \otimes \pi^W$ on $X^r$ induce maps $\pi^{V\otimes W}$ on $V\otimes W$ which give this space the structure of persistence module. We call this module \textit{the tensor product} of persistence modules $(V,\pi^V)$ and $(W,\pi^W)$. Another way to think of $V \otimes W$ is that $(V \otimes W)^r$ is the colimit in the category of (finite-dimensional, as is easy to see) vector spaces over our ground field of the diagram with objects $\{V^{s} \otimes W^{t}\}_{s + t \leq r}$ and maps $\pi_{s_1,s_2}\otimes \pi_{t_1,t_2} : V^{s_1} \otimes W^{t_1} \to V^{s_2} \otimes W^{t_2}$ for $s_1 \leq s_2$ and $t_1 \leq t_2$ (we use the convention that $\pi_{t,t}=\id_{V^t}$).

It is easy to see that we can also define the tensor product $f\otimes g:V\otimes W \rightarrow V'\otimes W'$ of persistence morphisms $f:V\rightarrow V'$ and $g:W \rightarrow W'$ by setting $f\otimes g([v_\alpha \otimes w_\beta])=[f(v_\alpha)\otimes g(w_\beta)]$.
\\
\\
Fixing a persistence module $W$ we get a functor $\otimes W:\pemod \rightarrow \pemod$ which acts on objects and morphisms by
$$\otimes W(V)=V\otimes W,~\otimes W(f)=f\otimes \id_W.$$
One can check that $\otimes W$ is a right exact functor and in order to define its derived functors we need to construct a projective resolution of every persistence module $V$. In the simplest case when $V$ is an interval module $Q((a,b])$ we have the following projective resolution of $V$ of length two:
$$0\rightarrow Q((b,+\infty)) \rightarrow Q((a,+\infty)) \rightarrow Q((a,b])\rightarrow 0,$$
where arrows denote obvious maps. Note that we used the fact that $Q((a,+\infty))$ is projective object for every $a\in \mathbb{R}$. One may also check that in fact a persistence module $V$ is projective if and only if its barcode contains no finite bars. Using this fact together with Theorem~\ref{StructureTheorem} we may construct a projective resolution of length two of every persistence module $V$ in the same manner as we did for the interval module. Recall that (classical) derived functors of $\otimes W$ applied to $V$ are computed as homologies of the sequence
$$\ldots \rightarrow P_2 \otimes W \xrightarrow{f_2 \otimes \id} P_1\otimes W \xrightarrow{f_1 \otimes \id} P_0 \otimes W \rightarrow 0,$$
where $\ldots \rightarrow P_2 \xrightarrow{f_2} P_1 \xrightarrow{f_1 } P_0 \xrightarrow{f_0} V \rightarrow 0$ is a projective resolution of $V$. Since every persistence module has a projective resolution of length two, there is only one non-trivial derived functor of $\otimes W$ which we denote by $Tor(\cdot, W)$. Both $\otimes$ and $Tor$ are symmetric in the sense that $V\otimes W \cong W \otimes V$ and $Tor(V,W)\cong Tor(W,V)$ and it immediately follows that if either $V$ or $W$ is projective $Tor(V,W)=0$.
\\
{\exm Let $V=Q((a,b]),W=Q((c,d])$ be two interval persistence modules. It follows directly from the definition of $\otimes$ that
$$V\otimes W=Q((a,b])\otimes Q((c,d])=Q((a+c, \min \{ a+d,b+c \}]).$$
In order to compute $Tor(Q((a,b]),Q((c,d]))$, let us take the following projective resolution of $Q((a,b])$:
$$0\rightarrow Q((b,+\infty)) \rightarrow Q((a,+\infty))\rightarrow Q((a,b]) \rightarrow 0.$$
After applying $\otimes Q((c,d])$ to this resolution we get
$$0\rightarrow Q((b+c,b+d]) \rightarrow Q((a+c,a+d])\rightarrow Q((a+c,\min \{ a+d,b+c \}]) \rightarrow 0,$$
and hence after calculating homology we get
$$Tor(Q((a,b]),Q((c,d]))=Q((\max \{ a+d,b+c \},b+d]).$$}
Our goal is to establish a K\"unneth type formula for filtered homology groups using $\otimes$ and $Tor$. Let us first recall the following definition.
{\defn We say that chain complex $(C_k,\partial_k), \partial_k:C_k\rightarrow C_{k-1},k\in \mathbb{Z}$ of finite dimensional vector spaces over a field $\mathbb{K}$ is \textit{filtered by function} $\nu$ if $\nu:C_* \rightarrow \mathbb{R}\cup \{ -\infty \}$ and
\begin{enumerate}
  \item[ 1)] $\nu(x)=-\infty$ if and only if $x=0$;
  \item[ 2)] For all $\lambda \in \mathbb{K},\lambda \neq 0$ it holds $\nu(\lambda x)=\nu(x)$;
  \item[ 3)] For all $x,y\in C_*$ it holds $\nu(x+y)\leq \max \{ \nu(x),\nu(y) \}$;
  \item[ 4)] For all $x\in C_*$ it holds $\nu (\partial_* x)\leq \nu(x)$.
\end{enumerate}
 }
{\rem This definition of chain complex filtered by function is the special case of the definition of {\it Floer-type complex} over Novikov field $\Lambda^{\mathbb{K},\Gamma}$ given in \cite{UsherZhang} in case $\Gamma=\{ 0 \}$ and valuation on $\mathbb{K}$ is trivial. \label{Remark:UZ}}
\\
\\
The main examples of filtered chain complexes of interest to us are Morse chain complex $CM_*(f)$ for Morse function $f$, where $f$ also serves as a filtration function and Floer chain complex $CF_*(H)_\alpha$ filtered by action functional $\mathcal{A}_H$, where $H$ is a Hamiltonian function and $\alpha$ is atoroidal or toroidally monotone class of free loops (see Section~\ref{Sect:Floer} for details).
\\
\\
Now if $(C_*,\partial_*,\nu)$ is a chain complex with filtration function $\nu$, we may define $C^t_*=\{x\in C_* | \nu(x)<t \}$ for evert $t\in \mathbb{R}$ and by property 4) we have that $\partial_*:C^t_*\rightarrow C^t_{*-1}$. This implies that $(C^t_*,\partial |_{C^t_*})$ is a new chain complex and we denote its homology by $H^t_*(C)$ and refer to it as {\it filtered homology}. Since $C^t\subset C^s$ for $t\leq s$ inclusions induce maps $\pi_{t,s}:H^t_*(C)\rightarrow H^s_*(C)$ which render $(H^t_*(C),\pi)$ into a persistence module. In order to obtain K\"unneth formula for filtered homology, we must examine the product of two filtered chain complexes. Let us start with an example.
\\
{\exm Let $(C^1_*,\partial^1_*,\nu^1)$ and $(C^2_*,\partial^2_*,\nu^2)$ be two filtered chain complexes given by
$$C^1_0=\langle x \rangle, C^1_1=\langle y \rangle, \partial^1 x=0, \partial^1 y=x, \nu^1(x)=a,\nu^1(y)=b,$$
and
$$C^2_0=\langle z \rangle, C^2_1=\langle w \rangle, \partial^1 z=0, \partial^1 w=z, \nu^2(z)=c,\nu^2(w)=d.$$
We have that
$$H^t_0(C^1)=Q^t((a,b]),H^t_0(C^2)=Q^t((c,d]),~H^t_*(C^1)=H^t_*(C^2)=0~\text{for}~*\neq 0.$$
The product complex $(C^1 \otimes C^2,\partial, \nu=\nu^1+\nu^2)$ (with usual product differential) is given by
$$(C^1 \otimes C^2)_0=\langle x \otimes z \rangle, (C^1 \otimes C^2)_1=\langle \{ x \otimes w, y\otimes z \} \rangle, (C^1 \otimes C^2)_2=\langle y \otimes w \rangle,$$
$$\partial (x\otimes z)=0, \partial (x\otimes w)=\partial (y\otimes z)=x\otimes z, \partial (y\otimes w)=x\otimes w - y\otimes z,$$
with filtration $\nu(x\otimes w)=a+d,\nu(y\otimes z)=b+c,\nu(x\otimes z)=a+b,\nu(y\otimes w)=b+d$. It readily follows that
$$ H_0^t(C^1\otimes C^2)=Q^t((a+b,\min \{ a+d,b+c \}])=(H_0(C^1)\otimes H_0(C^2))^t, $$
$$ H_1^t(C^1\otimes C^2)=Q^t((\max \{ a+d,b+c \},b+d])=(Tor(H_0(C^1), H_0(C^2))^t, $$
and $H_2^t=0$.
\label{Exm:Intervals}
}
\\
\\
Note that the $Tor$ functor naturally appears even in the simplest case of product of interval modules. As already mentioned, in this case torsion comes from finite bars in the barcode and hence is unavoidable even when we work with fields and vector spaces. We may now formulate the full statement.
{\prop[K\"unneth formula for filtered homology] Let $(C^1,\partial^1,\nu^1)$ and $(C^2,\partial^2,\nu^2)$ be two filtered chain complexes and let $(C^1\otimes C^2,\partial, \nu=\nu_1+\nu_2)$ be their product complex. Then for every $k\in \mathbb{Z}$ there exists a short exact sequence of persistence modules
$$0\rightarrow \bigoplus_{i+j=k} (H_i(C^1)\otimes H_j(C^2))^t \xrightarrow{K} H_k^t(C^1 \otimes C^2) \rightarrow$$
$$\rightarrow \bigoplus_{i+j=k-1} (Tor(H_i(C^1),H_j(C^2)))^t \rightarrow 0,$$
which splits, where $K$ denotes canonical map given by $K([\sum_i \lambda_i x_i]\otimes[\sum_j \mu_j y_j])=[\sum_{i,j} \lambda_i \mu_j x_i \otimes y_j]$.
\label{KunnethFormula}}
\\
\\
\textbf{Sketch of the proof.} We already saw in Example~\ref{Exm:Intervals} that the statement holds when $C^1$ and $C^2$ have the following form
$$\ldots \rightarrow 0 \rightarrow \langle y \rangle \rightarrow \langle x \rangle \rightarrow 0 \rightarrow \ldots $$
It readily follow that the statement is also true if we allow $C^1$ and $C^2$ to also be of the following form
$$\ldots \rightarrow 0 \rightarrow \langle x \rangle \rightarrow 0 \rightarrow \ldots$$
By Remark~\ref{Remark:UZ} we may look at our complexes as a special case of definition given in \cite{UsherZhang} and we may use the existence of singular value decomposition of operator $\partial$ proven there. This theorem essentially states that every filtered chain complex decomposes into direct sum of the simple complexes which have one of the two forms described above. Now, the general case follows from reduction to two simple ones and considerations about interval modules. \qed
\\
{\rem Essentially the same computation of the product of chain complexes as one presented in the Example~\ref{Exm:Intervals} and in the proof of Proposition~\ref{KunnethFormula} appears in \cite{Zhang}. The context is, however, slightly different, since we eventually work on the level of homology, while the author of \cite{Zhang} works on chain level. One may also try to prove Proposition~\ref{KunnethFormula} directly, without referring to much more general machinery developed in \cite{UsherZhang}.}

\subsection{Persistence modules with operators}\label{Sect:pmod-operator}

The methods we use, which are of independent interest, have to do with persistence modules endowed with an additional structure, and their equivariant version.

Consider the category $\pemodop$ with objects pairs $(V,A)$ with $V \in \pemod,$ and $A:V \to V[c_A],$ for certain $c_A \in \R,$ a morphism of persistence modules. Morphisms between $(V,A)$ and $(W,B),$ when $c_A = c_B$ consist of morphisms $F: V \to W$ of persistence modules, such that $F[c_A] \circ A = B \circ F,$ and if $c_A \neq c_B,$ only of the zero morpishm $V \to W.$

\subsubsection{Examples}

{\exm\textbf{(Shift operator)}
For each $\delta \geq 0,$ each $V \in \pemod$ comes with a canonical shift operator $sh(\delta): V \to V[\delta].$ For $\delta = 0,$ this is simply the identity operator. For $\delta >0,$ $sh(\delta)_t : V^t \to V^{t+\delta}$ is defined as the persistence structure map $\pi_{t,t+\delta}$ of $V.$ Hence $(V,sh(\delta))$ is an object of $\pemodop.$}

{\exm\textbf{($\mathbb{Z}_k$-action)}
Fix an integer $k \geq 2.$ Given a $\Z_k = \Z/k\Z$-representation in $\pemod,$ the action of the cyclic generator $1 \in \Z_k$ gives an operator $A: V \to V,$ with $c_A = 0$ (that satisfies $A^k= \id_V$).}

{\exm\textbf{(Product map)}
Consider a Morse function $f:X \to \R$ on a closed finite dimensional manifold $X$ of dimension $\dim X =m$. It defines a $\Z$-graded persistence module by $V_*(f)^s = H_*(\{f\leq s\},\K) = H_*(\{f < s\},\K),$ for $s$ a regular value of $f.$ Let $p_s: H_*(X) \to H_*(X,\{f \geq  s\}) =  H_*(\{f\leq s\},\{f = s\})$ be the natural map. Taking a class $a \in H_r(X),$ the intersection product with $p_s(a),$ $(p_s(a) \cap): V_*(f)^s \to V_{*+r-m}(f)^s$ defines an operator $(a \cap): V_*(f) \to V_*(f),$ with $c_{a \cap} = 0,$ that shifts the grading by $r-m.$\label{ProductExample}}

\subsubsection{Key estimate}

For two objects $(V,A)$ and $(W,B)$ of $\pemodop$ with $c_A = c_B,$  and $\delta \geq 0,$ define an operator-$\delta$-interleaving between them to be a $\delta$-interleaving $f: V \to W[\delta],$ $g: W \to V[\delta]$ that commutes with the operators $A$ and $B,$ that is \[f[c_A] \circ A = B[\delta] \circ f,\; g[c_B] \circ B = B[\delta] \circ g.\]

Define the operator-interleaving distance between them by \[d_{\opint}((V,A),(W,B)) = \inf \{\delta \geq 0|\, \text{there exists a } \delta\text{-operator-interleaving} \}.\]

\begin{prop}\label{prop:key-estimate}
For all $(V,A),$ $(W, B)$ in $\pemodop$ with $c_A =c_B,$ \[d_{inter}(\im(A),\im(B)) \leq d_{\opint}((V,A),(W,B)).\]
\end{prop}

Put $c:=c_A =c_B.$ The proof is an immediate diagram chase in the diagram (and its analogue with $f,g$ interchanged):

\begin{equation}\label{Diagram: Push gamma general and continuations}
  \displaystyle    \begin{array}{ccccc}
        V^t  & \xrightarrow{f} & W[\delta]^t & \xrightarrow{g[\delta]} & V[2\delta]^t\\
         \scriptstyle{A_t}{\downarrow} &   &  \scriptstyle{B_t[\delta]}{\downarrow} & & \scriptstyle{A_t[2\delta]}{\downarrow}\\
       V[c]^t   & \xrightarrow{f[c]} & W[c+\delta]^t & \xrightarrow{g[c+\delta]} & V[c+2\delta]^t\\
      \end{array}
   \end{equation}

\subsubsection{Discussion}

While Proposition \ref{prop:key-estimate} is elementary, it turns out to be useful already in the more basic examples.

{\exm\textbf{(Shift operator)}
Proposition \ref{prop:key-estimate} applied to the example of persistence shift maps, reduces to the following statement. If $V, W$ are $\delta$-interleaved, then $V'=\im sh(c)_V,\, W'=\im sh(c)_W$ are $\delta$-interleaved for every $c\in \mathbb{R}$. The reason is that with respect to shift operators, operator-$\delta$-interleaving is the same as $\delta$-interleaving, so $V, W$ being $\delta$-interleaved implies that they are also operator-$\delta$-interleaved.}

{\exm\textbf{(Intersection product)}
In Section \ref{Sect: Morse} we give examples of two Morse functions $f,g$ on a surface $\Sigma_2$ of genus $2$ with identical barcodes, and identical spectral invariants, the images of whose persistence modules under the intersection product with a class in $H_1(\Sigma_2,\mathbb{K})$ are, however, at a positive interleaving distance $c>0$. We conclude, by Proposition \ref{prop:key-estimate}, that any two functions in the respective orbits of $f,g$ under the indentity component of the diffeomorphism group are at $C^0$-distance $c>0.$ Indeed for such a diffeomorphism $\psi \in \Diff_0(\Sigma_2),$ $(V(f\circ \psi),a\cap)$ and $(V(f),a\cap)$ are isomorphic objects in $\pemodop,$ and still $d_{\opint}((V(f),a\cap),(V(g),a \cap)) \leq |f-g|_{C^0}.$ Indeed, the relevant interleavings commute with $a \cap.$ \label{IntersectionProduct}}


\subsection{Example of a Morse function on $\mathbb{T}^2\sharp\mathbb{T}^2$}\label{Sect: Morse}

We present an example in Morse homology which illustrates the effect of a product on Floer persistence module which we will define later and we also justify claims of Example~\ref{IntersectionProduct}.
\\
\\
Adopting the setup of Example~\ref{BasicExample} and Example~\ref{ProductExample}, to a Morse function $f$ on a closed manifold $X$ of dimension $m$ we associate a persistence module $(V^t_*(f), \pi)$ by taking $V^t_*(f)=H_*( \{ f<t \}, \mathbb{K} )$, the structure maps $\pi_{s,t}$ being induced by inclusion of sublevel sets. Alternatively, we may consider the Morse chain complex induced by critical points whose critical value is less than $t$. Now, $a\in H_*(M)$ acts on $V^t_*(M)$ by intersecting cycles (or by counting $\mathrm{Y}$-shaped configurations of gradient flow lines in Morse picture) and we get a map:
$$a\cap : V^t_*(f) \rightarrow V^t_{*+\deg a -m}(f).$$
Let $\Sigma_2$ be a surface of genus 2. We construct two Morse functions on $\Sigma_2$ which have same barcodes and same spectral invariants associated to every homology class, but their intersection barcodes with a fixed class differ by a finite bar. First, observe that $\Sigma_2 = \mathbb{T}^2\sharp\mathbb{T}^2$ and hence $H_1(\Sigma_2) \cong H_1(\mathbb{T}^2)\oplus H_1(\mathbb{T}^2)$, where generators are given by standard generators of $\mathbb{T}^2=S^1\times S^1$, namely two circles. We consider a Morse function $f:\Sigma_2\rightarrow \mathbb{R}$ given by the height function on the following picture:
\\
\begin{center}
\includegraphics[width=11cm]{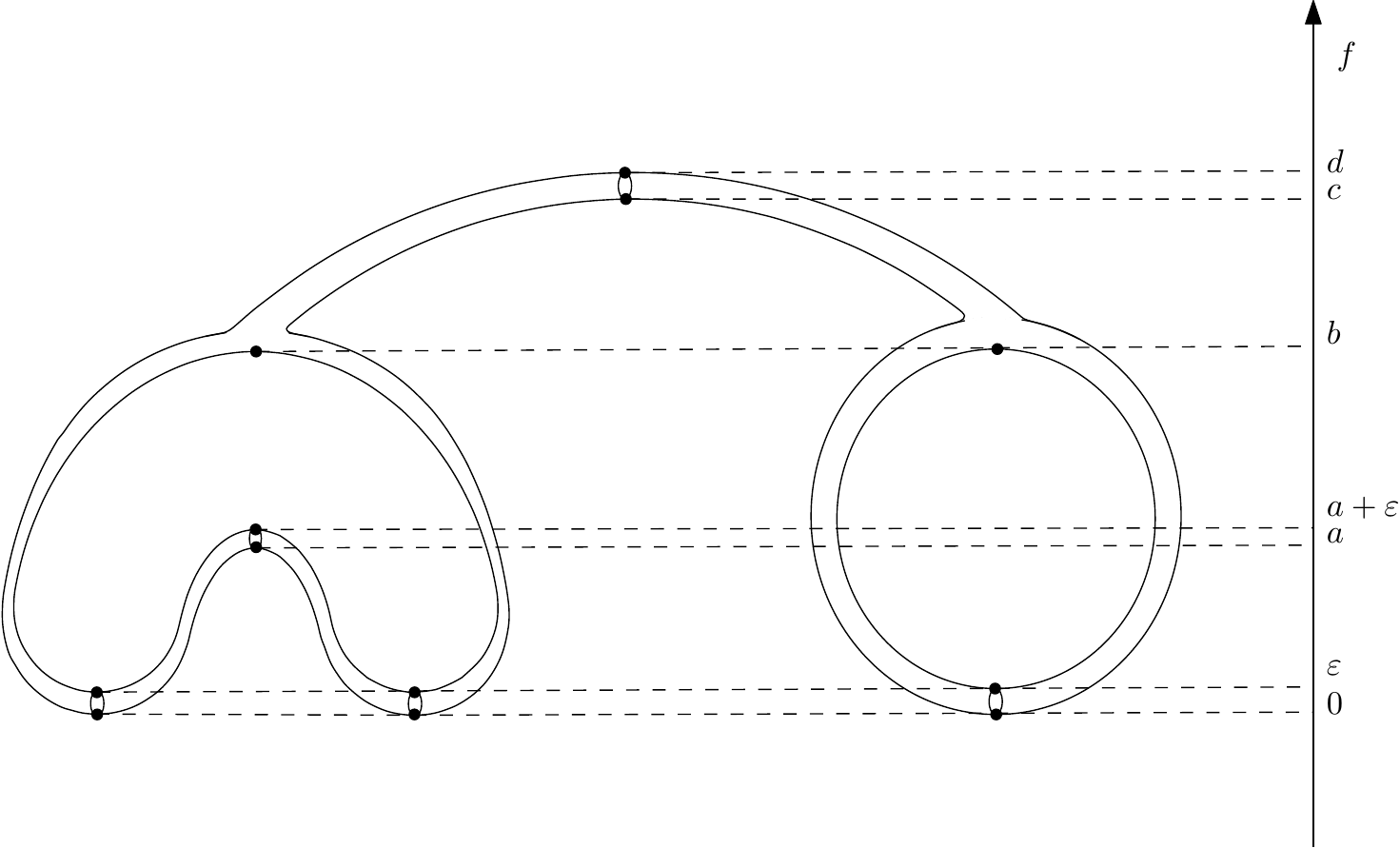}
\end{center}
We observe that $H_1(\Sigma_2)$ is generated by four homology classes represented by embedded circles, two of which have spectral invariants associated to $f$ equal to $\varepsilon$ and the other two with spectral invariants equal to $b$.
\\
\\
The other function we consider is the height function $g$ on the same picture with left and right reversed. More precisely, $g=f\circ \varphi$ where $\varphi: \Sigma_2 \to \Sigma_2$ is a diffeomorphism which interchanges two copies of $\mathbb{T}^2 \setminus D^2$ which we glue together to form $\Sigma_2$. Since $g=f\circ \varphi$, the barcodes of $f$ and $g$ are the same and they look as follows:
\begin{center}
\includegraphics[width=11cm]{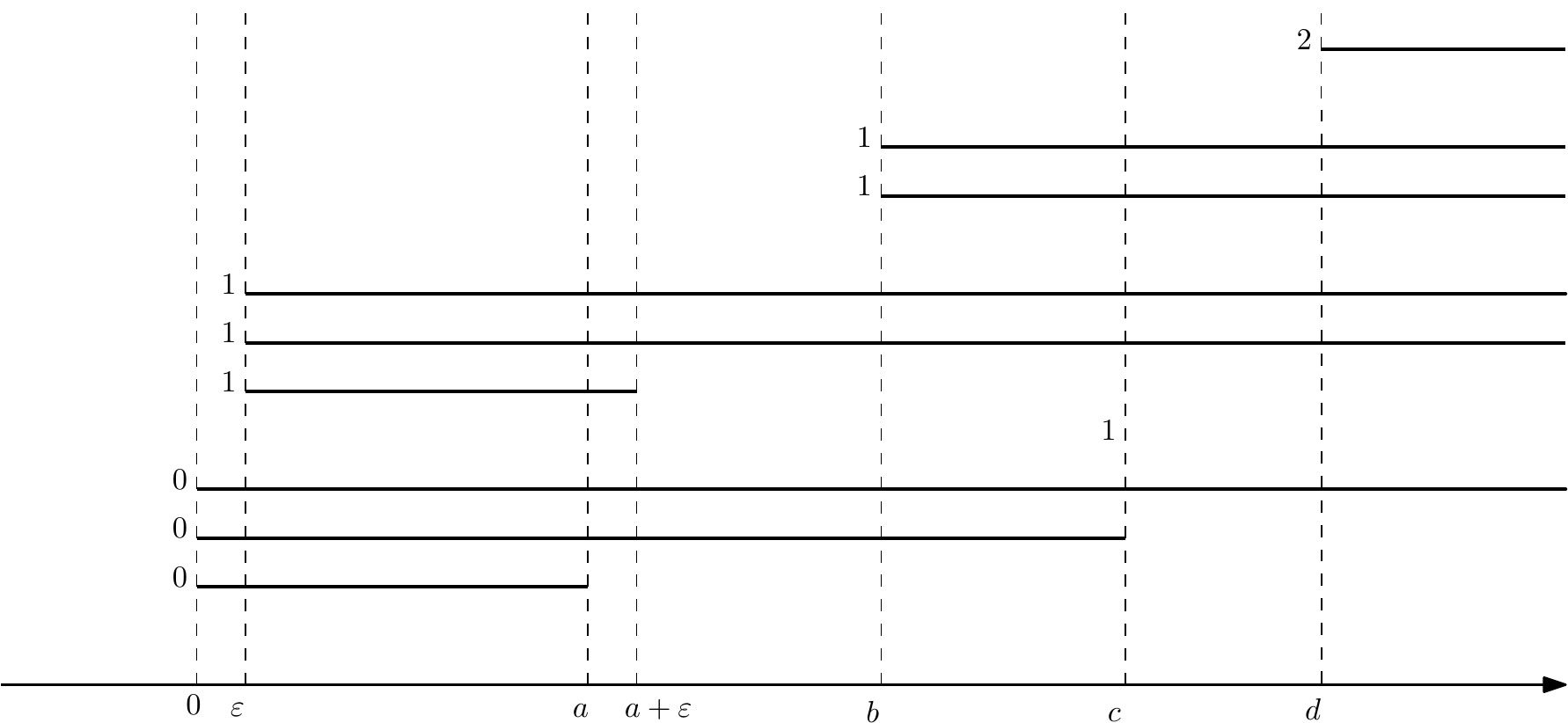}
\end{center}
One also readily checks that for every $z\in H_*(M)$, $c(z,f)=c(z,g)$ where $c(z,f),~c(z,g)$ are spectral invariants associated to functions $f$ and $g$ and a homology class $z$. This means that standard methods, namely barcodes and spectral invariants fail to distinguish between $f$ and $g$. However, after intersecting with one of the two big circles (for example the one on the left in the above picture), which corresponds to the homology class $e$ with spectral invariants $c(e,f)=c(e,g)=b$, we get the following intersection barcodes:
\begin{center}
\includegraphics[width=11cm]{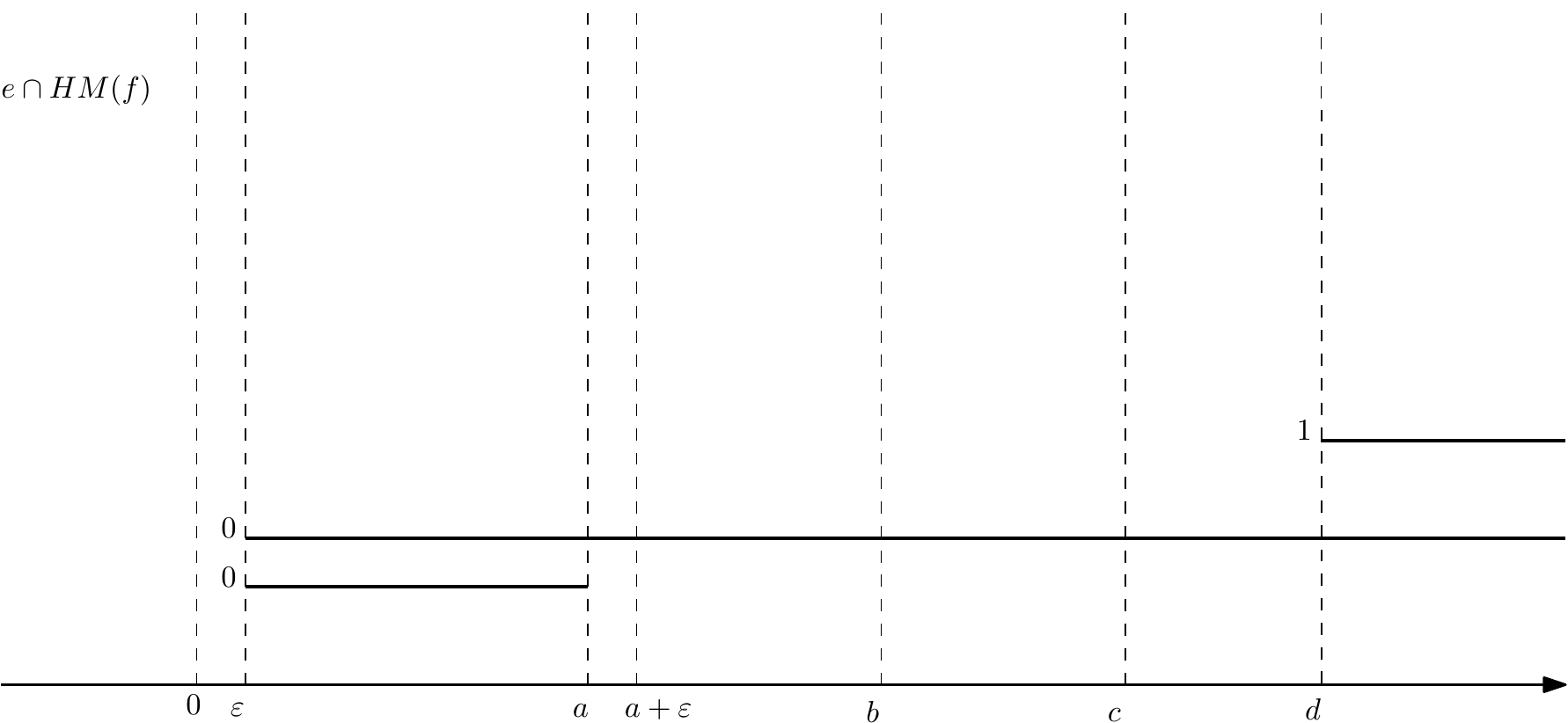}
\end{center}

\begin{center}
\includegraphics[width=11cm]{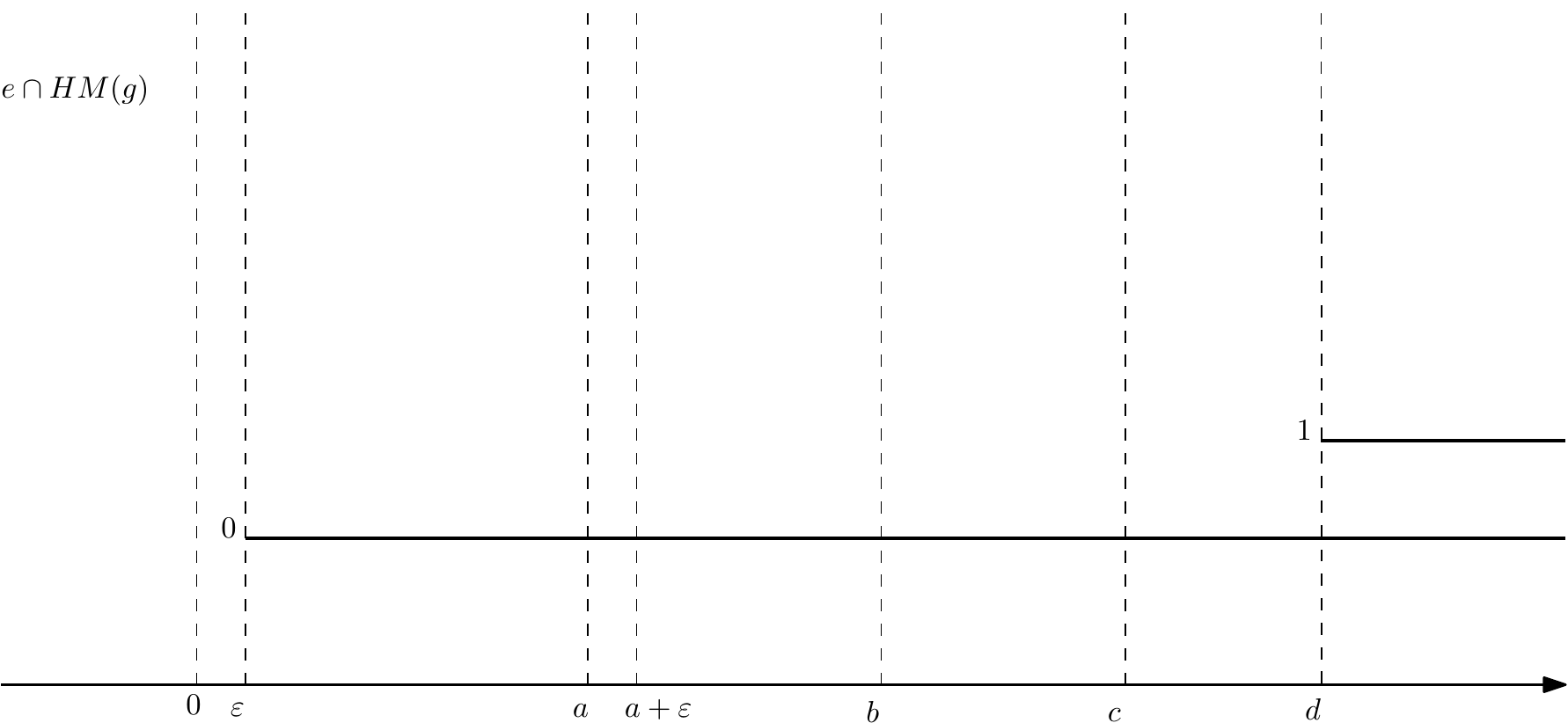}
\end{center}
These barcodes differ by a finite bar $(\varepsilon,a]$. Thus, by using the product structure in homology and analysing its effect on the barcode we are able to make a distinction between $f$ and $g$. Note also that the bar $(\varepsilon,a]$ did not exist in the original barcode.

It would be interesting to find a general formula for the image persistence module of the intersection by homology class $a.$ Examples show that this is not a trivial question.

\subsection{Equivariant version}\label{Sect:Equivariant}

In order to study the question of Hofer's distance to autonomous Hamiltonian diffeomorphisms and more generally to full $p$-th powers in $\Ham$, persistence modules with additional $\mathbb{Z}_p = \Z/p\Z$ action were used in \cite{PolShe}. A {\it $\mathbb{Z}_p$ persistence module } $(V,\pi,T)$ is a persistence module $(V,\pi)$ together with an automorphism $T:(V,\pi)\rightarrow (V,\pi)$ which satisfies $T^p=\id$. This definition immediately implies that $T_t:V^t \rightarrow V^t$ is a linear operator whose eigenvalues are $p$-th roots of unity. Hence, for $\zeta^p=1$, $\pi_{st}$ maps a $\zeta$-eigenspace of $T_s$ to $\zeta$-eigenspace of $T_t$ and we can define a $\zeta$-eigenspace of $T$ to be a persistence module obtained by restricting $\pi$ to $\zeta$-eigenspaces of each $T_t$.

We require the following immediate statement.

{\lem Let $(V_r,T_r)$, $r\in 1, \ldots , l$ be $\mathbb{Z}_p$ persistence modules, $(V,T)=(\bigoplus\limits_{r=1}^l V_r,\bigoplus\limits_{r=1}^l T_r)$ and denote by $L_\zeta$, $\zeta$-eigenspace of $T$, where $\zeta^p=1$, $\zeta \neq 1$. Then
$$L_\zeta = \bigoplus_{r=1}^l L_\zeta^r,$$
where $L_\zeta^r$ are $\zeta$-eigenspaces of $V_r$. \label{eigenspace}}

\medskip

Interleavings between $\mathbb{Z}_p$ persistence modules which commute with the $\mathbb{Z}_p$ action are called {\it equivariant}. Again, taking infimum over all $\delta>0$ such that $V$ and $W$ are eqivariantly $\delta$-interleaved gives us and {\it equivariant interleaving distance} which we denote by $\widehat{d}_{inter}(V,W)$. It immediately follows that
$$\widehat{d}_{inter}(V,W)\geq d_{inter}(V,W) ~\text{ and }~ \widehat{d}_{inter}(V,W) \geq d_{inter}(L_\zeta,K_\zeta),$$
where $L_\zeta$ and $K_\zeta$ are the $\zeta$-eigenspaces of $T_V$ and $T_W$ respectively.
\\
\\
Applying our new method to the equivariant situation is paramount to studying $\mathbb{Z}_p$ persistence modules with an operator $A:V\rightarrow V[c_A],$ which moreover, commutes with the $\Z_p$-action. Examples of such operators will come from a version of the pair-of-pants product in Floer homology.

{\defn A $\mathbb{Z}_p$ persistence module with an operator is a pair $(V,A)$ where $V$ is a $\mathbb{Z}_p$ persistence module and $A:V\rightarrow V[c_A]$ is a morphism of persistence modules that commutes with the $\mathbb{Z}_p$-action.}

\medskip
Let $(V,A)$ and $(W,B)$ be two $\mathbb{Z}_p$ persistence modules with operators with $c=c_A=c_B,$ and suppose that $f:V\rightarrow W[\delta]$ and $g:W\rightarrow V[\delta]$ is an equivariant $\delta$-interleaving. We say that this interleaving is {\it op-equivariant} if it respects the operator actions, that is
$$B(\delta)\circ f=f(c)\circ A,~~ A(\delta)\circ g=g(c)\circ B.$$
Taking infimum over all $\delta$ such that $V$ and $W$ are op-equivariantly $\delta$-interleaved gives us a new distance which we denote $\widehat{d}_{op}(V,W)$. Since $A$ and $B$ are $\mathbb{Z}_p$ persistence module morphisms we have that $\im(A)\subset V^a$ and $\im(B)\subset W^a$ are $\mathbb{Z}_p$ persistence submodules of $V^a$ and $W^a$. Every op-equivariant interleaving between $V$ and $W$ induces an equivariant interleaving between $\im(A)$ and $\im(B)$ which in particular implies
\begin{equation}\label{opinter_ineq}
\widehat{d}_{\opint}((V,A),(W,B))\geq \widehat{d}_{inter}(\im(A),\im(B)).
\end{equation}
Note however that in general this may not be an equality (see Section \ref{Sect: Morse}). 

{\rem The situation which we encounter when working with singular, Morse or Floer homology is not exactly the same as described above since our product map may change the degree and not just the filtration. One can overcome this ambiguity by giving a slightly more general definition analogous to the one given above, where $A:V^t \rightarrow \bar{V}^{t+a}$ for different persistence modules $V$ and $\bar{V}$ or by considering graded vector spaces.}
\\
\\
In order to tackle the problem of Hofer's distance to full powers in $\Ham$ a numerical invariant $\mu_p(W)$ called {\it multiplicity sensitive spread} was defined in  \cite{PolShe}. We recall the definitions and properties of $\mu_p$ and an auxiliary invariant $\mu_{p,\zeta}$ which we use later (see  \cite{PolShe} for proofs).
\\
\\
Let $\cB$ be a barcode, $I$ and interval and denote by $m(\cB,I)$ the number of bars in $\cB$ containing $I$ (counted with multiplicities).
We will write $\mu_p(\cB)$ for a supremum of those $c \geq 0$ for which there exists an interval $I$ of length greater than $4c$ such that $m(\cB, I) = m(\cB, I^{2c})=l$ with $l$ not divisible by $p$. Using this notation we define $\mu_{p,\zeta}$ as
$$\mu_{p,\zeta}(W)=\mu_p(\mathcal{B}(L_\zeta)),$$
where $L_\zeta$ is $\zeta$-eigenspace of $T$. Now $\mu_p$ is defined as
$$\mu_p(W)=\max_{\zeta} \mu_{p,\zeta}(W).$$
We have that
$$\vert \mu_p(\mathcal{B}(L_\zeta)) - \mu_p(\mathcal{B}(K_\zeta)) \vert = \vert \mu_{p,\zeta}(V) - \mu_{p, \zeta}(W) \vert \leq d_{bottle}(\mathcal{B}(L_\zeta),\mathcal{B}(K_\zeta)),$$
where $L_\zeta$ and $K_\zeta$ are the $\zeta$-eigenspaces of $T_V$ and $T_W$ respectively.
\\
\\
A $\mathbb{Z}_p$ persistence module $(W,T)$ is called {\it a full p-th power} if $T=S^p$ for some morphism $S:W \rightarrow W$.

\medskip

From now on we impose the same assumption on the ground field $\K$ as in Section~\ref{Sect:Hofer}. An important property of $\mu_p$ for such a ground field $\K$ is that $\mu_p(W)=0$ given that $W$ is a full $p$-th power.

\section{Floer theory and Hofer's geometry}

\subsection{Product map on Floer persistence module}\label{Sect:Floer}

Let $(M,\omega)$ be a closed symplectic manifold, and denote by $c_1(TM)$ the first Chern class of the tangent bundle, equipped with any $\omega$-compatible almost complex structure. Take a homotopy class of free loops $\alpha \in \pi_0(\mathcal{L}M)$ and denote by $\mathcal{L}_\alpha M$ all loops in class $\alpha$. We say that $(M,\omega)$ is $\alpha$-toroidally monotone if there exists $\kappa>0$ such that
$$\langle [\omega], A \rangle = \kappa \cdot \langle c_1(TM), A \rangle, $$
for all $A\in Im (\Psi)$, where $\Psi : \pi_1(\mathcal{L}_\alpha M) \rightarrow H_2(M, \mathbb{Z})$ sends a loop $\beta \in \pi_1( \mathcal{L}_\alpha M )$, regarded as a map $\beta : \mathbb{T}^2 \rightarrow M$, to $\beta_* ([\mathbb{T}^2])$. It readily follows that $M$ is also spherically monotone with same monotonicity constant $\kappa$, that is
$$ [\omega]=\kappa \cdot c_1(TM),$$
where both $[\omega]$ and $c_1(TM)$ are regarded as functionals on $\pi_2(M)$. Assuming $M$ is $\alpha$-toroidally monotone\footnote{All the considerations in this section also apply to $\alpha$ which is symplectically atoroidal, meaning $\omega=c_1=0$ on $\pi_1(\mathcal{L}_\alpha M)$.}, to every element $\tilde{f}\in \widetilde{\Ham}(M)$ of the universal cover of $\Ham(M)$, that is non-degenerate in class $\al,$ we associate a Floer persistence module $HF_*^t(\tilde{f})_\alpha$ with parameter $t$ (see \cite{PolShe,UsherBD2}). Taking $\alpha = [pt]$, and $\tilde{f}=\id \in \widetilde{\Ham}(M)$ we recover, for fixed degree and large $t,$ the usual Floer homology for monotone manifolds, which is canonically isomorphic to $QH(M)$ (see \cite{PSS}). We can define a product map
\begin{equation}
\label{eq:product} *:HF^t_*(\id)_{pt} \otimes HF^s_*(\tilde{f})_\alpha \rightarrow HF^{t+s}_*(\tilde{f})_\alpha,\end{equation}
by counting pairs of pants on the chain level.
\\
\\

We remark that here, and later on in Section \ref{Sect:stab}, we deal with degenerate Hamiltonian diffeomorphisms $\id \in \Ham(M)$ and $\phi \times \id \in \Ham(\Sigma \times N)$ (for $\Sigma$ a closed symplectic surface of higher genus, and $N$ a monotone symplectic manifold) which, however, are of Morse-Bott degeneracy in the appropriate classes of orbits. Associating a persistence module to this situation can be treated in a number of ways. First, arguing up to epsilon everywhere, we could replace $\id$ by the flow of a sufficiently $C^2$-small Morse function $h,$ considered as a Hamiltonian (see e.g. \cite[proof of Prop. 4.2]{PolShe} for arguments of this type). In this new, perturbed, setting a persistence module is defined. Moreover, we could fix $h,$ replace $\id$ with the flow of $\delta \cdot h,$ and look at the appropriate persistence modules as $\delta \to 0.$ It is easy to see, by use of action estimates in PSS maps \cite{PSS}, for example, that they converge in interleaving distance to a well-defined genuine persistence module (which is uniquely determined up to isomorphism by this property). In other words, \cite[Definition 2.8]{PolShe} applies in this case, and gives a persistence module in the sense described above, since the set of the critical values of the action functional of the zero Hamiltonian is discrete, in our case. Finally, one could use Frauenfelder's approach of cascades \cite{Frauenfelder} to the Morse-Bott case, which readily yields a persistence module by the same procedure as in \cite{PolShe, UsherBD2}. In the case of $\id \in \Ham(M),$ the last two approaches compute $HF^t_r(\id)_{pt}$ in degree $r \in \Z$ as follows. The Novikov field $\Lambda_{\K}$ admits a non-Archimedean valuation $\nu: \Lambda_{\K} \to \R\cup \{-\infty\},\; \sum\limits_{n\in \mathbb{Z}}a_nq^n \mapsto \max \{ n \cdot (\kappa c_N)|\; a_n \neq 0\}.$ This valuation naturally extends to $QH(N),$ by declaring that $\nu(x) = 0$ for all non-zero $x \in H_*(N,\K) \otimes 1.$ Then $HF^t_r(\id)_{pt} = QH(M)_r^t$ that is defined as $QH^t_r(N)=\{x\in QH_r(N) |\, \nu(x)<t \}.$

In the first two approaches, the action estimates for the product map follow from \cite[Section 4.1]{Schwarz:action-spectrum},\cite{EntovInvent},\cite[Section 6.2]{Oh-construction}. In the third approach, the product map takes the form of counting "spiked cylinders", quite similar to the definition of the PSS map \cite{PSS} (see e.g. \cite{Char} and references therein for details on the more complicated, Lagrangian, version).

Let us examine some of the properties of this product.
\\
\\
Denote by $\tilde{d}$ the Hofer pseudo-distance on $\widetilde{\Ham}(M)$ and by $d$ the Hofer distance on $\Ham(M)$. We write $\tilde{f}\in \widetilde{\Ham}(M)$ for a homotopy class of paths relative endpoints in $\Ham(M)$ and $f\in \Ham(M)$ for its endpoint. Let $\nu: QH(M) \rightarrow \mathbb{R}$ be the natural valuation. Since $HF_*(\id)_{pt} \cong QH_*(M)$, fixing homogeneous $a\in QH(M)$ we obtain a map
$$a*:HF^t_r(\tilde{f})_\alpha \rightarrow HF^{t+\nu(a)}_{r-2n+\deg a}(\tilde{f})_\alpha.$$
The map $a\ast$ is a persistence module morphism between $V^t_r=HF^t_r(\tilde{f})_\alpha$ and $\widetilde{V}^t_r=HF^{t+\nu(a)}_{r-2n+\deg a}(\tilde{f})_\alpha$. Moreover, it follows from standard considerations in Floer theory that $a*$ commutes with continuation maps
$$C(F,G):HF^t_r(F)_\alpha \rightarrow HF^{t+\mathcal{E}^+(G-F)}_r(G)_\alpha,$$
where $\mathcal{E}^+(G-F)=\int_{0}^{1}\max_M(G_t-F_t)dt$.
\\
\\
Now, let $g\in \Ham(M)$ and define a map
$$P(g):HF^t_*(\tilde{f})_\alpha \rightarrow HF^t_*(\tilde{g\circ f\circ g^{-1}})_\alpha,$$
by acting with $g$ on all the objects appearing in the construction of Floer chain complex. More precisely, on the chain level $P(g)$ defines an isomorphism of filtered chain complexes
$$P(g):(CF(H,J)_\alpha, \mathcal{A}_H) \rightarrow ((CF(H\circ g^{-1}, g_*(J))_\alpha, \mathcal{A}_{H\circ g^{-1}})),$$
by sending a periodic orbit $z(t)$ of $H$ to a periodic orbit $g(z(t))$ of $H\circ g^{-1}$. This map is called {\it the push-forward map} (see \cite{PolShe} for a detailed treatment of push-forward maps). One can check that $P(g)$ and $a*$ commute.
\\
\\
Our objects of interest are Floer persistence modules of the form $HF^t_*(\tilde{f}^p)_\alpha$ for $\tilde{f}\in \widetilde{\Ham}(M)$. In this case $P(f):HF^t_*(\tilde{f}^p)_\alpha \rightarrow HF^t_*(\tilde{f}^p)_\alpha$ defines a $\mathbb{Z}_p$ action on $HF^t_*(\tilde{f}^p)_\alpha$ and we get a $\mathbb{Z}_p$ Floer persistence module. Since $P(f)$ and $a*$ commute, $a*$ is a $\mathbb{Z}_p$ persistence module morphism and we wish to treat it as an operator on $HF^t_*(\tilde{f}^p)_\alpha$ and apply considerations from Section \ref{Sect:Equivariant}. To do so, define a $\mathbb{Z}_p$ persistence module
$$W^t_r(a,\tilde{f}^p)=\im(a*)=(a*)(HF^t_r(\tilde{f}^p)_\alpha) \subset HF^{t+\nu(a)}_{r-2n+\deg a}(\tilde{f}^p)_\alpha,$$
with $\mathbb{Z}_p$ action given by $P(f)$. Denote by $F_t$ and $G_t$ normalized 1-periodic Hamiltonians generating paths in $\Ham(M)$ which represent classes of $\tilde{f}$ and $\tilde{g}$ in $\widetilde{\Ham}(M)$ respectively and by $F^{(p)}_t=pF_{pt}$ and $G^{(p)}_t=pG_{pt}$ normalized 1-periodic Hamiltonians generating paths which represent $\tilde{f}^p$ and $\tilde{g}^p$. Continuation maps
$$HF^t_r(F^{(p)})_\alpha \xrightarrow{C(F^{(p)},G^{(p)})} HF^{t+p\cdot \mathcal{E}^+(G-F)}_r(G^{(p)})_\alpha,$$
and
$$HF^{t+p\cdot \mathcal{E}^+(G-F)}_r(G^{(p)})_\alpha \xrightarrow{C(G^{(p)},F^{(p)})} HF^{t+p\cdot(\mathcal{E}^+(G-F)-\mathcal{E}^-(G-F))}_r(F^{(p)})_\alpha,$$
induce a $p\cdot(\mathcal{E}^+(G-F)-\mathcal{E}^-(G-F))$ op-equivariant interleaving between $HF^t_*(\tilde{f}^p)_\alpha$ and $HF^t_*(\tilde{g}^p)_\alpha$, where $\mathcal{E}^-(G-F)=\int_{0}^{1}\min_M(G_t-F_t)dt$. Taking infimum over all $F$ and $G$ generating $\tilde{f},\tilde{g}\in \widetilde{\Ham}(M)$ we get that
$$p\cdot \tilde{d}(\tilde{f},\tilde{g})\geq \widehat{d}_{\opint}(HF^t_*(\tilde{f}^p),HF^t_*(\tilde{g}^p)),$$
which together with (\ref{opinter_ineq}) gives us
\begin{equation}\label{Hofer_Estimate}
p\cdot \tilde{d}(\tilde{f},\tilde{g})\geq \widehat{d}_{\opint}(HF^t_*(\tilde{f}^p),HF^t_*(\tilde{g}^p)) \geq \widehat{d}_{inter}(W^t_*(a,\tilde{f}^p),W^t_*(a,\tilde{g}^p)).
\end{equation}

{\rem Let $f\in \Ham(M)$ and fix a lift $\tilde{f}\in \widetilde{\Ham}(M)$ of $f$. We can use $W^t_r(a,\tilde{f}^p)$ to estimate Hofer's distance from $f$ to $p$-th powers inside $\Ham(M)$. Indeed, denote by $\Powers_p(M) \subset \Ham(M)$ the set of all $p$-th powers of Hamiltonian diffeomorphisms and by $\widetilde{\Powers_p}(M) \subset \widetilde{\Ham}(M)$ the set of all lifts of elements from $\Powers_p(M)$. In other words $\widetilde{\Powers_p}(M) = \pi^{-1}(\Powers_p(M))$ under the natural projection $\pi:\widetilde{\Ham}(M) \to \Ham(M).$ For $\tilde{g}\in \widetilde{\Powers_p}(M)$, we have that $W^t_r(a,\tilde{g}^p)$ is a full $p$-th power persistence module because $g=\phi^p$ implies $P(\phi)^p=P(g)$ and $P(\phi)$ restricts to $W^t_r(a,\tilde{g}^p)$ because $P(\phi)$ and $a*$ commute. It follows that $\mu_p(W^t_r(a,\tilde{g}^p))=0$ and thus
$$\vert \mu_p(W^t_r(a,\tilde{f}^p)) \vert = \vert \mu_p(W^t_r(a,\tilde{f}^p)) - \mu_p(W^t_r(a,\tilde{g}^p)) \vert \leq \widehat{d}_{inter}(W^t_r(a,\tilde{f}^p),W^t_r(a,\tilde{g}^p)),$$
which together with (\ref{Hofer_Estimate}) gives us
$$\vert \mu_p(W^t_r(a,\tilde{f}^p)) \vert \leq \widehat{d}_{inter}(W^t_r(a,\tilde{f}^p),W^t_r(a,\tilde{g}^p)) \leq p\cdot \tilde{d}(\tilde{f},\tilde{g}).$$
Finally, we have
\begin{equation}
d(f,\Powers_p(M))=\tilde{d}(\tilde{f},\widetilde{\Powers_p}(M))\geq \frac{1}{p} \cdot \vert \mu_p(W^t_r(a,\tilde{f}^p)) \vert.
\end{equation}
\label{distance}}

\subsection{Stabilization and the egg-beater example}\label{Sect:stab}

We now turn to a manifold $M$ of the form $M=\Sigma \times N$, where $\Sigma$ is surface of genus at least 4 and $N$ is spherically monotone symplectic manifold with monotonicity constant $\kappa$. The element $\tilde{\psi}_\lambda\in \widetilde{\Ham}(M)$ which we consider is
$$\tilde{\psi}_\lambda=\tilde{\varphi}^p_\lambda \times \id,~\tilde{\varphi}_\lambda\in \widetilde{\Ham}(\Sigma),~\id \in \widetilde{\Ham}(N),$$
where $\tilde{\varphi}_\lambda$ is given by the {\it egg-beater flow} on $\Sigma$, with mixing parameter $\lambda$. Construction and detailed analysis of egg-beater flow are carried out in \cite{Team,PolShe}. What we will use is that there exists a family of Hamiltonian flows $\tilde{\varphi}_\lambda$ on $\Sigma$, depending on an unbounded increasing real parameter $\lambda$, along with a family of classes of free loops $\alpha_\lambda$ on $\Sigma$ which satisfy:
\begin{enumerate}
  \item[1)] $\varphi_\lambda^p$ has exactly $2^{2p}$ $p$-tuples of fixed points with same indices and actions $\{ z,\varphi_\lambda(z),\ldots, \varphi^{p-1}_\lambda(z) \}$, for each $\lambda$;
  \item[2)] If $z_1$ and $z_2$ belong to different $p$-tuples their action differences satisfy $$|\mathcal{A}(z_1)-\mathcal{A}(z_2)|\geq c_0 \lambda + O(1);$$
  \item[3)] The indices of all fixed points are bounded by a constant which does not depend on $\lambda$.
\end{enumerate}
The class $\overline{\alpha}_\lambda\in \pi_0(\mathcal{L}M)$ which we consider is a product of classes
$$\overline{\alpha}_\lambda = \alpha_\lambda \times pt, ~\alpha_\lambda \in \pi_0(\mathcal{L}\Sigma), $$
$\Sigma$ being symplectically $\alpha_\lambda$-atoroidal. Our manifold $M$ will be $\overline{\alpha}_\lambda$-toroidally monotone with same monotonicity constant $\kappa$. We will leave out these classes from the notation and write $HF^t_*(\tilde{\varphi}_\lambda^p \times \id)$ and $HF^t_*(\tilde{\varphi}_\lambda^p)$ for $HF^t_*(\tilde{\varphi}_\lambda^p \times \id)_{\overline{\alpha}_\lambda}$ and $HF^t_*(\tilde{\varphi}_\lambda^p)_{\alpha_\lambda}$. Let us now work out the example which proves Theorem~\ref{main}.
{\prop Let $\tilde{\varphi}_\lambda$ be the egg-beater flow and assume $e\in QH(N)$ satisfies assumptions of Theorem~\ref{main}. There exists $k\in \mathbb{Z}$ such that
$$\mu_p(W^t_k([\Sigma]\otimes e, \tilde{\varphi}_\lambda^p \times \id)) \geq c\lambda + O(1),$$
for some $c>0$, when $\lambda \rightarrow +\infty$. Here $[\Sigma]\otimes e \in QH(M)=QH(\Sigma \times N).$ \label{MainProp}}
\\
\\
{\proof Let $\alpha_1,\alpha_2$ be two toroidally monotone classes of free loops in symplectic manifolds $M_1$ and $M_2$, with same monotonicity constant $\kappa$ (we may also take one of both of them to be atoroidal) and let $\tilde{\phi}\in \widetilde{\Ham}(M_1),\tilde{\psi}\in \widetilde{\Ham}(M_2)$. The manifold $M_1 \times M_2$ is symplectic and the class $\alpha_1 \times \alpha_2$ is toroidally monotone with the same monotonicity constant $\kappa$. Now, we apply Proposition~\ref{KunnethFormula} for general filtered homologies to Floer chain complexes filtered by action functional and Floer persistence modules to get the short exact sequence:
$$
0\rightarrow \bigoplus_{i+j=k} (HF_i(\tilde{\phi})_{\alpha_1}\otimes HF_j(\tilde{\psi})_{\alpha_2})^t \xrightarrow{K} HF_k^t(\tilde{\phi} \times \tilde{\psi})_{\alpha_1\times \alpha_2} \rightarrow$$
$$\rightarrow \bigoplus_{i+j=k-1} (Tor(HF_i(\tilde{\phi})_{\alpha_1},HF_j(\tilde{\psi})_{\alpha_1}))^t \rightarrow 0,$$
for $K([\sum_i \lambda_i x_i]\otimes[\sum_j \mu_j y_j])=[\sum_{i,j} \lambda_i \mu_j x_i \otimes y_j]$.
\\

In our case $\tilde{\phi}=\tilde{\varphi}_\lambda^p$, $\alpha_1=\alpha_\lambda$, $\tilde{\psi}=\id\in \widetilde{\Ham}(N),\alpha_2=\{pt\}$ and we have $HF^t_*(\id)_{\{pt\}}=QH^t_*(N)$, where $QH^t_*(N)=\{x\in QH_*(N) |\, \nu(x)<t \}$ is a persistence module with trivial structure maps given by $\pi_{s,t}(x)=x$ since $QH^s_*(N)\subset QH^t_*(N)$ for $s\leq t$.

This readily gives us that the barcode of $QH^t_*(N)$ has only infinite bars and thus $QH^t_*(N)$ is a projective persistence module and $Tor(HF_i(\tilde{\varphi}_\lambda^p),QH_j(N))=0$ for all $j\in \mathbb{Z}$, which implies that
$$K:\bigoplus_{i+j=k} (HF_i(\tilde{\varphi}_\lambda^p)\otimes QH_j(N))^t \rightarrow HF_k^t(\tilde{\varphi}_\lambda^p \times \id),$$
is an isomorphism. Moreover, it holds that $P(\varphi_\lambda \times \id)\circ K=K\circ (P(\varphi_\lambda)\otimes \id)$ (see \cite{PolShe} for a proof in the atoroidal case, the proof in the toroidally monotone case is the same) and thus $K$ is also an isomorphism of $\mathbb{Z}_p$ persistence modules. Now, consider multiplication by $e$ as a persistence module morphism $(e*):QH^t_r(N)\rightarrow QH^{t+\nu(e)}_{r-2n+\deg e}(N)$ between $QH^t_r(N)$ and shifted module $QH^{t+\nu(e)}_{r-2n+\deg e}(N)=QH^t_{r-2n+\deg e}(N)[\nu(e)]$, for every $r \in \mathbb{Z}$. Our product map splits on the components of the product, i.e. it enters the following commutative diagram:
\begin{center}
\begin{tikzcd}
\bigoplus\limits_{i+j=k} (HF_i(\tilde{\varphi}_\lambda^p)\otimes QH_j(N))^t \arrow[r,"K"] \arrow[d, "\id\otimes (e*)"]
& HF_k^t(\tilde{\varphi}_\lambda^p \times \id) \arrow[d,"( \lbrack\Sigma \rbrack \otimes e)*"] \\
\bigoplus\limits_{i+j=k} (HF_i(\tilde{\varphi}_\lambda^p)\otimes QH_{j-2n+\deg e}(N)[\nu(e)])^t \arrow[r,"K"]
& HF^t_{k-2n+\deg e}(\tilde{\varphi}_\lambda^p \times \id)[\nu(e)]
\end{tikzcd}
\end{center}
where each arrow represents a $\mathbb{Z}_p$ persistence module morphism. Using this diagram we calculate
$$W^t_k=W^t_k([\Sigma]\otimes e, \tilde{\varphi}_\lambda^p \times \id)=\bigoplus_{r\in I_k}(HF_{k-r}(\tilde{\varphi}_\lambda^p)\otimes (e*)(QH_r(N)))^t,$$
where $I_k$ is the set of all $r$ such that there exists a fixed point of $\tilde{\varphi}_\lambda^p$ of index $k-r$ and $\im(e*)^t=(e*)(QH_r^t(N))\subset QH^t_{r-2n+\deg e}(N)[\nu(e)].$ Let us describe the barcode of $\im(e*)^t$ explicitly. 
\\
\\
First, note that we have inclusions of all $QH^t_r(N)$ into full quantum homology $QH_r(N)=QH^{+\infty}_r(N)$ and moreover for $t\leq s$, $QH_r^t(N)\subset QH_r^s(N)\subset QH_r(N)$ and structure maps act as $\id$ under these inclusion. Now, $E_r=e*(QH_r(N))\subset QH_{r-2n+\deg e}$ is the image of full quantum homology group $QH_r(N)$ and by the assumption $\dim_{\mathbb{K}}E_r=b_r(e)$.
\\
\\
We may also look at $E_r$ as a persistence submodule of the shifted module $E_r^t\subset QH^t_{r-2n+\deg e}[\nu(e)]$ and $(e*):QH^t_r(N)\rightarrow E_r^t$ is a persistence module morphism. Since structure maps on $E_r^t$ are restrictions of structure maps on $QH_{r-2n+\deg e}(N)[\nu(e)]$, we again have that they act as $\id$ under the inclusions to full quantum homology group $QH_{r-2n+\deg e}(N)$ and the same holds for $\im(e*)^t$. This implies that the barcode of $\im(e*)^t$ contains no finite bars. Now, if we denote $a_r=\min \{\nu(x)|x\in QH_r(N) \}$ and $A_r=\max \{\nu(x)|x\in QH_r(N) \}$, it follows that $(e*)(QH^t_r(N))=0$ for $t\leq a_r$ and $(e*)(QH^t_r(N))=E_r$ for $t>A_r$ and thus the barcode of $\im(e*)^t$ consists of bars $(c_{r,1},+\infty),\ldots,(c_{r,b_r(e)},+\infty)$ where $a_r\leq c_{r,1} \leq \ldots \leq c_{r,b_r(e)} \leq A_r$. Moreover, since $\mathbb{Z}_p$ action on $QH_r^t$ is trivial for all $r$ we have that
$$\im(e*)^t = \bigoplus_{i=1}^{b_r(e)} (Q^t((c_{r,i},+\infty)),\id),$$
as $\mathbb{Z}_p$ persistence modules, which together with the above diagram gives us
$$(W^t_k, P(\varphi_\lambda \times \id))\cong \bigoplus_{r\in I_k} \bigg(\bigg(HF_{k-r}(\tilde{\varphi}_\lambda^p) \otimes \bigoplus_{i=1}^{b_r(e)}Q((c_{r,i},+\infty))\bigg)^t, P(\varphi_\lambda)\otimes \id \bigg) ,$$
isomorphism being given by $K$. Elementary calculations on interval persistence modules now imply
$$(W^t_k, P(\varphi_\lambda \times \id)) \cong \bigoplus_{r\in I_k} \bigoplus_{i=1}^{b_r(e)} \bigg( HF^{t-c_{r,i}}_{k-r}(\varphi^p_\lambda)_{\alpha_\lambda}, P(\varphi_\lambda) \bigg). $$
Denoting the $\zeta$ eigenspace of $( HF^t_k(\varphi^p_\lambda)_{\alpha_\lambda}, P(\varphi_\lambda))$ by $L^t_{k,\zeta}$ and $\zeta$ eigenspace of $(W^t_k,P(\varphi_\lambda \times \id))$ by $L^t_\zeta$ we have by Lemma~\ref{eigenspace}
$$L^t_\zeta \cong \bigoplus_{r\in I_k} \bigoplus_{i=1}^{b_r(e)} L^{t-c_{r,i}}_{k-r,\zeta}.$$
The indices of fixed points of egg-beater map are uniformly bounded (the bound does not depend on $\lambda$) and thus we have $|r|<M$ for $r\in I_k$ for some constant $M$ not depending on $\lambda$. This also gives us that there exist a constatn $C>0$ independent of $\lambda$ such that $|a_r|<C$ and $|A_r|<C$ for all $r\in I_k$ and thus $|c_{r,i}|<C$ for all $r\in I_k,i=1,\ldots, b_r(e)$. By Lemma~\ref{sumpersistence} we have that
$$d_{bottle} \bigg(\mathcal{B}(L^t_\zeta) , \mathcal{B} \bigg( \bigoplus_{r\in I_k} ( L^t_{k-r,\zeta})^{b_r(e)} \bigg) \bigg) <C,$$
and hence by Lipschitz property of $\mu_p$ we have
$$\mu_p(W^t_k)\geq \mu_{p,\zeta}(W^t_k)=\mu_p \bigg( \bigoplus_{r\in I_k} \bigoplus_{i=1}^{b_r(e)} L^{t-c_{r,i}}_{k-r,\zeta} \bigg) \geq \mu_p \bigg( \bigoplus_{r\in I_k} ( L^t_{k-r,\zeta})^{b_r(e)} \bigg) -C.$$
Assume now that $p \nmid b_{r_0}(e) $ and that the index of a fixed point $z_0$ of $\varphi^p_\lambda$ with minimal action $A=A(z_0)$ in class $\alpha_\lambda$ is $d_0$. Taking $k=d_0+r_0$ we have that
$$ \bigoplus_{r\in I_k} (L^t_{k-r,\zeta})^{b_r(e)} = (L^t_{d_0,\zeta})^{b_{r_0}(e)} \oplus \bigoplus_{r \neq r_0} (L^t_{k-r,\zeta})^{b_r(e)}.$$
If $z$ is a fixed point of $\varphi^p_\lambda$ with action $A(z)\neq A$ it follows that $A(z)\geq B = A+c_0 \lambda + O(1)$ and we have that
$$m \bigg( \mathcal{B} \bigg( \bigoplus_{r\in I_k} (L^t_{k-r,\zeta})^{b_r(e)} \bigg) , (A,B] \bigg)=b_{r_0(e)}.$$
Now, $p \nmid b_{r_0}(e)$ and thus
$$\mu_p \bigg( \mathcal{B} \bigg( \bigoplus_{r\in I_k} (L^t_{k-r,\zeta})^{b_r(e)} \bigg) \bigg) \geq \frac{c_0}{4}\lambda + O(1)$$
which gives us $\mu_p(W^t_k)\geq c\lambda + O(1)$ as claimed. \qed
}
\\
\\
The proof of Theorem~\ref{main} follows directly from Proposition~\ref{MainProp} and Remark~\ref{distance}.

\subsection{Erratum: behavior $\mu_p$ under stabilization in the aspherical case.}

This erratum is written in order to correct a mistake in Theorem 4.24 in \cite{PolShe}. The main theorem (which this mistake could potentially affect), \cite[Theorem 1.3]{PolShe}, holds still. See Theorem E1 and the update to the proof of \cite[Theorem 1.3]{PolShe} below.

Alternatively, as noted in Example \ref{exa2:main-thm}, \cite[Theorem 1.3]{PolShe} holds as a special case of the main theorem, Theorem \ref{main} of the current paper, and its proof extends the proof of \cite[Theorem 1.3]{PolShe}.

In fact the estimate $\mu_p(\phi) \leq \mu_p(\phi \times \id)$ cannot be expected to hold, as can be seen by elementary examples. The error in the proof of Theorem 4.24 is contained in the implication "Thus we are left with $i=0$..." because the barcodes $\cB_{r-i}(\phi)_\zeta$ for $i > 0$ can have $I$ and $I^{2c}$ with different multiplicities, thus affecting the value of $\mu_{p,\zeta}(r,\phi \times \id_N).$

Denote  \[\gamma_{p,\zeta}(r,\phi) = \frac{1}{2} \max_{i>0} \beta(\cB_{r-i}(\phi)_\zeta).\]
By (26), and the remarks on the Kunneth formula in the proof of Theorem 4.24, it is immediate that \[\mu_{p,\zeta}(r,\phi \times \id_N) \geq \mu_{p,\zeta}(r,\phi) - \gamma_{p,\zeta}(r,\phi).\] Indeed \[d_{bottle}(\cB_{r}(\phi)_\zeta,\cB_{r}(\phi \times \id_N)_\zeta) \leq \gamma_{p,\zeta}(r,\phi),\] which can be seen by erasing all intervals corresponding to ($b_i(N)$-copies of) the barcode $\cB_{r-i}(\phi)_\zeta$ (recall that $\beta(\cB)$ is the maximal length of a finite bar in the barcode $\cB$).

Thus denoting \[\mu_{p,\zeta}^{\mathrm{reduced}}(r,\phi) =  \mu_{p,\zeta}(r,\phi) - \gamma_{p,\zeta}(r,\phi),\] and \[\mu_{p}^{\mathrm{reduced}}(r,\phi) = \max_{\zeta}\mu_{p,\zeta}^{\mathrm{reduced}}(r,\phi)\] we replace Theorem 4.24 by the following.

\begin{thme}\label{thm--mu2-stabilization-correction}

For $\phi \in \Ham(M),$ $\al \in \pi_0(\cL M),$ and any closed connected symplectically aspherical manifold $N,$ consider the stabilization $\phi \times \id \in {\Ham(M \times N)}$ of $\phi.$ Then we have \[\mu_p^{\mathrm{reduced}}(\phi) \leq \mu_p(\phi \times \id_N) \leq \mu_p(\phi),\] the value $\mu_p(\phi \times \id_N)$ being computed in the class $\al \times pt_N$ in $\pi_0(\cL (M \times N)).$

\end{thme}

Now we turn to Section 5.1 and show how to adapt the proof of Theorem 1.3 in view of the corrected Theorem E1 above. The necessary changes are:

\begin{itemize}
\item The sentence

\medskip

"Further, among the $2^{2p}$ $p$-tuples of fixed points of $\phi_\lambda^p$ in the class $\alpha_\la$
choose the $p$-tuple, say $\{z,\phi_\lambda(z),\ldots, (\phi_\la)^{p-1}(z)\}$ with the minimal action. Let $r$ be the index of $z.$"
\medskip

should be corrected to

\medskip
"Further, among the $2^{2p}$ $p$-tuples of fixed points of $\phi_\lambda^p$ in the class $\alpha_\la$
choose the $p$-tuple, say $\{z,\phi_\lambda(z),\ldots, (\phi_\la)^{p-1}(z)\}$ with the minimal index $r,$ and minimal action among $p$-tuples of this index."
\medskip

\item The passage

\medskip

"By the definition of the multiplicity-sensitive spread, we conclude that $\mu_p(\phi_\lambda) \geq \la (c-2\varepsilon )/4$"

\medskip
should read
\medskip

 "By the definition of the multiplicity-sensitive spread and the observation that $\gamma_{p,\zeta}(r,\phi_\lambda) = 0,$ we conclude that $\mu_p^{\mathrm{reduced}}(\phi_\lambda) \geq \la (c-2\varepsilon )/4$"

 \medskip
\end{itemize}

\section*{Acknowledgements}

We thank Iosif Polterovich for critical remarks on the manuscript, and Jelena Kati\'c, Jovana Nikoli\'c, Darko Milinkovi\'c, Igor Uljarevi\'c, and Jun Zhang for useful discussions.

Work on this paper was carried out while E.S. was staying at the Institute for Advanced Study. He thanks the IAS, and Helmut Hofer in particular, for their warm hospitality.

\bibliographystyle{abbrv}

{\footnotesize
\bibliography{Stabilization_Paper}}

\bigskip

\noindent
\begin{tabular}{lll}
Leonid Polterovich & Egor Shelukhin & Vuka\v sin Stojisavljevi\'c\\
School of Mathematical Sciences & IAS, Princeton, and & School of Mathematical Sciences\\
Tel Aviv University & DMS at U. of Montreal  & Tel Aviv University\\
polterov@post.tau.ac.il & egorshel@gmail.com & vukasin@post.tau.ac.il \\
\end{tabular}

\end{document}